\documentclass[hidelinks,final]{siamart220329}

\usepackage{smalljeff}
\usepackage{booktabs,makecell}
\usepackage{stmaryrd}
\usepgfplotslibrary{fillbetween, groupplots}
\usepgfplotslibrary{colorbrewer}
\usepgfplotslibrary{external}
\usetikzlibrary{decorations.markings}
\usepackage{datatool}
\usepackage{xfrac}

\usepackage{afterpage,placeins}

\newcommand{\E}{\mathbb{E}}
\DeclareMathOperator{\Cov}{Cov}

\tikzset{->-/.style={decoration={
  markings,
  mark=at position #1 with {\arrow{>}}},postaction={decorate}}
}

\newcommand\undermat[2]{%
  \makebox[0pt][l]{$\smash{\underbrace{\phantom{%
    \begin{matrix}#2\end{matrix}}}_{\text{$#1$}}}$}#2}

\newcommand\writemat[2]{%
\pgfplotstableread[header= false, col sep=comma]{#1}\mytable
\pgfplotstabletypeset[
  begin table = {},
  end table = {},
  skip coltypes,
  write to macro=\mymatrix,
  assume math mode = true,
  typeset=false,
  every head row/.style={output empty row},
  every column/.style={fixed zerofill={true}, precision=6},
  #2
]{\mytable}
\begin{bmatrix}
\mymatrix
\end{bmatrix}
}

\newcommand\sidds{{\scriptscriptstyle \spadesuit}}
\newcommand\lsoi{{\scriptscriptstyle \clubsuit}}
\newcommand\generic{{\scriptscriptstyle \diamondsuit}}

\newcommand{\TheTitle}{Simultaneous Identification and Denoising of Dynamical Systems}
\newcommand{\TheAuthors}{Jeffrey M. Hokanson, Gianluca Iaccarino, and Alireza Doostan}
\headers{Simultaneous Identification and Denoising}{\TheAuthors}

\title{{\TheTitle}\thanks{Submitted to the editors DATE.
\funding{This material is based upon work supported by the Department of Energy, National Nuclear Security Administration under Award Number DE-NA0003968.}}
}

\author{Jeffrey M. Hokanson\thanks{
	Department of Aerospace Engineering Sciences, University of Colorado at Boulder, Boulder, CO 80309 
	(jeffrey@hokanson.us, alireza.doostan@colorado.edu)}
\and Gianluca Iaccarino\thanks{
Center for Turbulence Research, Stanford University, Stanford, CA 94305 (jops@stanford.edu)
}
\and Alireza Doostan\footnotemark[2]
}

\begin{document}
\maketitle
\begin{abstract}
In recent years there has been a push to discover the governing equations dynamical systems
directly from measurements of the state,
often motivated by systems that are too complex to directly model.
Although there has been substantial work put into
such a discovery, 
doing so in the case of large noise has proved challenging.
Here we develop an algorithm for
the Simultaneous Identification and Denoising of a Dynamical System (SIDDS).
We infer the noise in the state measurements
by requiring that the denoised state
satisfies the dynamical system with an equality constraint.
This contrasts to existing work 
where the mismatch in the dynamics is added as a penalty in the objective.
Assuming the nonlinear differential equation is represented in a pre-defined basis,
we develop sequential quadratic programming approach to solve the SIDDS problem
featuring a direct solution of KKT system with a specialized preconditioner.
We also show how to add a sparsity promotion regularization into SIDDS
using an iteratively reweighted least squares approach.
Our resulting algorithm obtains estimates of the dynamical system 
that achieve the Cram\'er-Rao lower bound up to discretization error.
This enables SIDDS to provide substantial improvements compared to existing techniques:
SIDDS substantially decreases the data burden for accurate identification,
recovers optimal estimates with lower sample rates,
and the sparsity promoting variant 
discovers the correct sparsity pattern with larger noise.
\end{abstract}
\begin{keywords}
	dynamical systems,
	model discovery,
	inverse problems,
	parameter estimation,
	sparse recovery
\end{keywords}
\begin{AMS}
	34A55, 
	65L09, 
	90C55, 
	93B30 
\end{AMS}
\begin{DOI}
\end{DOI}
\section{Introduction}
We consider the problem of identifying a dynamical system
from measurements of its state.
Suppose the state of this system $\ve x \in\R^d$ satisfies a first order, 
autonomous ordinary differential equation (ODE)
\begin{equation}\label{eq:ode}
	\left\lbrace
	\begin{aligned}
		\dve x(t) &= \ve f(\ve x(t)), \\
		\ve x(0) &= \ve x_0, 
	\end{aligned}
	\right. \quad \text{where} \quad \ve f: \R^d\to \R^d.
\end{equation}
Our goal is to identify $\ve f$
given access to $m$ noisy observations of the state
$\lbrace \ve y_j \rbrace_{j=1}^m$,
where $\ve y_j \approx \ve x(t_j)$
at times $\lbrace t_j \rbrace_{j=1}^m \subset [0, T]$.
This problem emerges in a variety of contexts
from model reduction~\cite{PW16} to system identification~\cite{BPK16}.

\subsection{Parameterization}
An important choice for recovering the operator $\ve f$
is its parameterization.
Here, we express $\ve f$ as a sum of $n$ scalar-valued basis functions $\phi_k :\R^d \to \R$
with corresponding coefficients $\ve c_k \in \R^d$ following~\cite{BPK16}
\begin{equation}\label{eq:expansion}
	\ve f(\ve x; \ma C ) \coloneqq \sum_{k=1}^n \ve c_k \phi_k(\ve x),
	\quad \text{where} \quad \ma C \coloneqq \begin{bmatrix} \ve c_1^\trans \\ \vdots \\ \ve c_n^\trans \end{bmatrix}
	\in \R^{n \times d}.
\end{equation}
A typical choice for $\lbrace \phi_k\rbrace_{k=1}^n$ is a polynomial basis;
e.g., total degree-2 basis in two dimensions ($d=2$) is:
\begin{align*}
	\phi_1(\ve x) &= 1,  & \phi_2(\ve x) &= x_1, & \phi_3(\ve x) &= x_2, & \phi_4(\ve x) &= x_1^2, &
	\phi_5(\ve x) &= x_1x_2,  & \phi_6(\ve x) &= x_2^2.
\end{align*}
Non-polynomial terms such as $\phi_k(\ve x) = \sin(x_1)$ can also be incorporated.
In some cases, we will further seek a parsimonious, interpretable expression for $\ve f$
by seeking a sparse coefficient matrix $\ma C$.
There are other parameterizations possible for $\ve f$,
for example, a neural network~\cite{RKB19}.

\subsection{A Naive Least Squares Approach}
A variety of methods exploit a simple linear relationship to estimate
the coefficients $\ma C$.
Suppose we have access to exact measurements of the state $\ve x(t_j)$
and its derivative $\dve x(t_j)$
and stack these into matrices
\begin{equation}
	\ma X \coloneqq \begin{bmatrix} \ve x(t_1)^\trans \\ \vdots \\ \ve x(t_m)^\trans \end{bmatrix} \in \R^{m \times d} 
	\quad \text{and} \quad 
	\dma X \coloneqq \begin{bmatrix} \dve x(t_1)^\trans \\ \vdots \\ \dve x(t_m)^\trans \end{bmatrix} \in \R^{m \times d}.
\end{equation}
Further, we build the matrix-valued function $\ma \Phi: \R^{m \times d} \to \R^{m\times n}$
containing the evaluations of the basis functions:
\begin{equation}
	\ma \Phi(\ma X) \coloneqq \begin{bmatrix}
		\phi_1(\ve x(t_1)) & \cdots & \phi_n(\ve x(t_1))\\
		\vdots & & \vdots \\
		\phi_1(\ve x(t_{m})) & \cdots & \phi_n(\ve x(t_m))
	\end{bmatrix} \in \R^{m\times n}.
\end{equation}
With this notation, based on the differential equation~\cref{eq:ode} and the expansion in~\cref{eq:expansion},
the coefficients $\ma C$ must satisfy the linear system:
\begin{equation}
		\dma X 
		= \begin{bmatrix} \dve x(t_1)^\trans \\ \vdots \\ \dve x(t_m)^\trans \end{bmatrix}
		= \begin{bmatrix} \ve f(\ve x(t_1))^\trans \\ \vdots \\ \ve f(\ve x(t_m))^\trans \end{bmatrix}
		= \begin{bmatrix} 
			\sum_{k=1}^n \ve c_k^\trans \phi_k(\ve x(t_1)) \\
			\vdots \\
			\sum_{k=1}^n \ve c_k^\trans \phi_k(\ve x(t_m))
		\end{bmatrix}
		= \ma \Phi(\ma X) \ma C.
\end{equation}
When this problem is well-posed, 
there is a unique solution for $\ma C$ given $\dma X$ and $\ma \Phi(\ma X)$;
numerically, $\ma C$ can be identified by solving the least squares problem
\begin{equation}\label{eq:exact}
	\min_{\ma C \in \R^{n \times d}} \| \dma X - \ma \Phi(\ma X) \ma C\|_\fro^2,
\end{equation}
where $\| \cdot \|_\fro$ denotes the Frobenius norm.
However, this problem may not have a unique solution;
for example, if there is a linear combination of basis functions 
encoding a conservation law for the system,
then $\ma \Phi(\ma X)$ will have a nontrivial nullspace.

The limitations of experimental measurements present two difficulties:
we may not have access to the derivative $\dve x(t_j)$
and our measurements $\ve y_j$ of $\ve x(t_j)$ are invariably contaminated by noise.
We can correct the former difficulty by estimating the derivative 
using finite-difference approximations based on measurements $\ve y_j$.
When measurements are uniformly spaced in time with time-step $\delta$,
e.g., $t_j = \delta(j-1)$,
we can build a finite difference matrix $\ma D\in \R^{m\times m}$
using a $q$-point central difference rule in the interior
and an order $q-1$ accurate rule on the boundary.
For example, a $3$-point finite difference matrix is
\begin{align}\label{eq:der3}
	\ma D =
		\frac{1}{2\delta}
		{\small
		\renewcommand{\arraystretch}{0.9}
		\begin{bmatrix}
			-3 & 4 & -1 \\
			-1 & 0 & \phantom{-}1 \\
			& \ddots & \ddots & \ddots \\
			&& \ddots & \ddots & \ddots \\
			&& & -1 & 0 & 1 \\
			&& & \phantom{-}1 & -4 & 3 
		\end{bmatrix}
		}
	\in \R^{m \times m}.
\end{align}
This provides the derivative approximation $\dma X \approx \ma D \ma X$.
To deal with the other difficulty of only having access to noisy measurements $\ve y_j$ of $\ve x(t_j)$, 
we replace $\ma X$ with $\ma Y$:
\begin{equation}
	\ma X \approx \ma Y \coloneqq \begin{bmatrix} \ve y_1^\trans \\ \vdots \\ \ve y_m^\trans \end{bmatrix} \in \R^{m \times d}
\end{equation}
Thus making the substitutions $\dma X \approx \ma D \ma X\approx \ma D \ma Y$
and $\ma \Phi(\ma X) \approx \ma \Phi(\ma Y)$, 
we can infer the coefficients $\ma C$
by solving the least squares problem
\begin{equation}\label{eq:lsopinf}
	\min_{\ma C \in \R^{n\times d}} \| \ma D \ma Y  - \ma \Phi(\ma Y)\ma C\|_\fro^2.
\end{equation}
We refer to this as \emph{Least Squares Operator Inference} (LSOI)
and this problem appears as a component in many approaches.

Although LSOI is inexpensive to solve, 
it inherits many limitations from the chain of approximations used in its construction.
The noise in $\ma Y$
may be amplified through the nonlinear basis functions in $\ma \Phi$
and the derivative estimate $\ma D \ma Y$ may be similarly inaccurate due to noise.
Thus when LSOI appears, 
it is often coupled with techniques to ameliorate these issues.
For example \emph{Sparse Identification of Nonlinear Dynamics} (SINDy)~\cite{BPK16}
introduces the desire for a sparse coefficient matrix $\ma C$.
This sparsity promotion improves the conditioning of the linear system
by (effectively) deleting columns from $\ma \Phi(\ma Y)$~\cite[Cor.~7.3.6]{HJ12}.
Another corrective measure is to find a better estimate of $\ma X$
by applying a denoising technique to the data $\ma Y$
to obtain an estimate $\hma Y$ that replaces $\ma Y$ in LSOI~\cref{eq:lsopinf}.
When $\hma Y$ is a better estimate of $\ma X$ than the original data $\ma Y$,
this improves estimate of $\ma C$;
however careful tuning is needed to avoid either
over smoothing or insufficiently removing noise~\cite{CPD22x}.
The only place where LSOI appears without embellishment
is in model reduction~\cite[eq.~(12)]{PW16}
where the data $\lbrace \ve y_j\rbrace_{j=1}^m$ comes directly from numerical simulations.

\subsection{Identification and Denoising}
A fundamental limitation of LSOI is that, 
in the presence of noise, 
the coefficient estimates do not satisfy the discretized dynamics;
that is, $\ma D \ma Y \ne \ma \Phi(\ma Y) \ma C$.
Our key contribution is to introduce an auxiliary variable $\ma Z$ 
which \emph{will} satisfy the discretized dynamics exactly: $\ma D \ma Z = \ma \Phi(\ma Z)\ma C$.
Then we minimize the mismatch between $\ma Z$ and the data $\ma Y$ in an appropriate norm.
Assuming noise follows independent and identically distributed (i.d.d.) normal distribution where
$\ve y_j \sim \mathcal{N}(\ve x(t_j), \sigma^2 \ma I)$,
we solve the constrained optimization problem
\begin{equation}\label{eq:sidds}
	\begin{split}
		\min_{\substack{\ma C\in \R^{n\times d}, \ \ma Z\in \R^{m\times d}}} & \ \| \ma Y - \ma Z\|_\fro^2 \\
		\text{s.t.}  \ \  \qquad & \ \ma D \ma Z = \ma \Phi(\ma Z) \ma C.
	\end{split}
\end{equation}
We call this \emph{Simultaneous Identification and Denoising of Dynamical Systems} (SIDDS).

Although we have introduced SIDDS using a similar framework to LSOI,
SIDDS can alternatively be derived from a more traditional inverse problem approach; see, e.g.~\cite{EAO00}.
SIDDS approximately solves the ODE-constrained optimization problem,
\begin{equation}\label{eq:ode_con}
	\begin{split}
	\min_{\ma C\in \R^{n\times d}, \ \ve z_0 \in \R^d, \ve \zeta : [0,T] \to \R^d} & \  \sum_{j=1}^m \| \ve y_j - \ve \zeta(t_j)\|_2^2  \\
	\text{s.t.} & \ \dve \zeta(t) = \ve f(\ve \zeta(t); \ma C), \quad \ve \zeta(0) = \ve z_0.
	\end{split}
\end{equation}
In particular, SIDDS~\cref{eq:sidds} uses the \emph{discretize-then-optimize} approach
to discretize \cref{eq:ode_con} in time
and uses a \emph{full-space} optimization approach
where both the coefficients $\ma C$ and predicted state history $\ve z_j = \ve \zeta (t_j)$
are variables in the optimization problem; see, e.g.~\cite{HA01}.
The full-space approach contrasts to a \emph{reduced-space} approach
where the predicted state history $\ve z_j$ is implicitly defined by $\ma C$ and $\ve z_0$; 
i.e., solving
\begin{equation}\label{eq:reduced}
	\min_{\ma C \in \R^{n \times d}, \ \ve z_0 \in \R^d}  \ 
		\sum_{j=1}^m \| \ve y_j - \ve \zeta(t_j; \ma C, \ve z_0)\|_2^2,
	\quad\text{where}\quad 
	\left\lbrace
	\begin{aligned}
		\dve \zeta(t; \ma C, \ve z_0) &= \ve f(\ve \zeta(t); \ma C), \\
		\ve \zeta(0; \ma C, \ve z_0) &= \ve z_0.
	\end{aligned}
	\right.
\end{equation}
The reduced-space problem is then solved using an unconstrained optimization algorithm
without storing the state history $\lbrace \ve z_j\rbrace_{j=1}^m$~\cite[Alg.~2.2.]{Hei13x}.
Although inexpensive, this approach can yield a hypersensitive objective
since when approximating a chaotic system,
small changes in $\ma C$ and $\ve z_0$ yield exponentially increasing changes 
in $\ve \zeta(t; \ma C, \ve z_0)$ with $t$.
This, in turn, implies that local information, such as gradients and Hessians, 
are not accurate beyond some small neighborhood
and makes finding a meaningful descent direction challenging.
In contrast, our full-space approach has a convex objective
and although the equality constraint needs to be satisfied at termination,
intermediate iterations do not need to satisfy the constraint exactly
allowing more effective exploration of the parameter space. 

The SIDDS problem shares similarities with problems in data assimilation~\cite{ABN16}.
Like data assimilation, we seek to estimate the true state from noisy measurements,
and as in some formulations, we seek to estimate the parameters of the underlying differential equation.
Unlike SIDDS which is solved using a full space method,
most data assimilation algorithms use a reduced space approach due to scaling concerns.

\subsection{A First Example}
Why use SIDDS instead of the existing approaches based on LSOI?
In short, SIDDS recovers far more accurate dynamical systems from noisy measurements.
As a first example, 
consider the simple harmonic oscillator $\ddot x(t) = -x(t)$ in first order form:
\begin{equation}\label{eq:sho}
	\begin{bmatrix}
		\dot{x}_1 \\ \dot{x}_2 
	\end{bmatrix}
	= 
	\begin{bmatrix}
	0 & 1 \\
	-1 & 0 
	\end{bmatrix}
	\begin{bmatrix}
		x_1 \\ x_2
	\end{bmatrix}
	\quad \text{with} \quad 
	\begin{bmatrix}
		x_1(0) \\ x_2(0)
	\end{bmatrix}
	=\begin{bmatrix} 1 \\ 0 \end{bmatrix}.
\end{equation}
Reconstructing this system,
we use two linear basis functions:
$\phi_1(\ve x) = x_1$ and $\phi_2(\ve x) = x_2$.
Then if we take $m=2000$ measurements with sample rate $\delta =10^{-2}$
that are contaminated with i.i.d.\ standard normal noise with unit covariance,
using the 3-point derivative approximation in~\cref{eq:der3}, we recover the systems
\begin{center}
\vspace{0.5em}
\noindent
\begin{tabular}{@{}cc}
	LSOI & SIDDS \\[0.5em]
	{\small $ 
	\begin{bmatrix}
		\dot{x}_1 \\ \dot{x}_2 
	\end{bmatrix}
	= 
	\begin{bmatrix}
	-0.23580696 & 2.85537791  \\
	 -2.60758784 &  0.2062968 
	\end{bmatrix}
	\begin{bmatrix}
		x_1 \\ x_2
	\end{bmatrix}
	$}
	&
	{\small $ 
	\begin{bmatrix}
		\dot{x}_1 \\ \dot{x}_2
	\end{bmatrix}
	= 
	\begin{bmatrix}
		-0.01155883 &  0.99851078 \\
		-1.0097105 &  -0.0019787 
	\end{bmatrix}
	\begin{bmatrix}
		x_1 \\ x_2
	\end{bmatrix}
	$}.
\end{tabular}
\vspace{0.5em}
\end{center}
Here SIDDS's estimate is accurate to 2-significant figures;
LSOI captures none!

\begin{figure}
\centering
\begin{tikzpicture}
\begin{groupplot}[
	group style = {group size = 1 by 2,
		horizontal sep = 2em,
		vertical sep = 1em,
		},
	width = .97\linewidth,
	height= 0.2\linewidth,
	xmin = 0, xmax = 20,
	ymin = -4, ymax = 4,
	clip mode=individual,
]
\pgfplotstableread{data/fig_sho.dat}\sho

\nextgroupplot[ylabel = {$x_1$}, xticklabels = {,,}]
\addplot[grey, mark=*, only marks, mark size = 0.2pt] table [x=t, y=y1] {\sho};

\addplot[black, line width=3pt] table [x=t, y=x1] {\sho};
\addplot[yellow, thick] table [x=t, y=z1] {\sho};

\nextgroupplot[ylabel = {$x_2$}, xlabel=time $t$]
\addplot[grey, mark=*, only marks, mark size = 0.2pt] table [x=t, y=y2] {\sho};

\addplot[black, line width=3pt] table [x=t, y=x2] {\sho};
\addplot[yellow, thick] table [x=t, y=z2] {\sho};

\end{groupplot}
\end{tikzpicture}

\caption{%
SIDDS accurately denoises the state and recovers the dynamical system
in the simple harmonic oscillator example.
Here the true system evolution is shown as a black line, 
the measured data as gray points,
and the recovered system as a yellow line.
}
\label{fig:sho}
\end{figure}
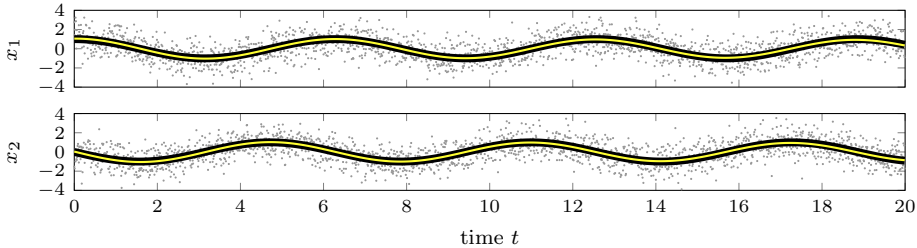

The reason why SIDDS performs so much better than LSOI is that,
subject to optimization algorithm finding the global minimizer,
$\ma Z$ is a maximum likelihood estimate of the true data $\ma X$
so that $\ma Z\approx \ma X$ as illustrated in \cref{fig:sho}.
With this accurate estimate of $\ma X$ in $\ma Z$,
SIDDS provides an accurate estimate of the coefficients $\ma C$.
In contrast, LSOI does not provide a maximum likelihood estimate of $\ma C$
due to the way it correlates noise in the expression $\ma D \ma Y - \ma \Phi(\ma Y) \ma C$. 
Details concerning these results are provided in \cref{sec:crb}.

\subsection{Sparsity Promotion\label{sec:intro:sparse}}
With an eye to obtain an easily interpretable 
expansion of $\ve f$ in the basis $\lbrace \phi_k \rbrace_{k=1}^n$ as in~\cref{eq:expansion},
it is common to apply a regularization to encourage the coefficients $\ma C$ to be sparse.
\emph{Sparse Identification of Nonlinear Dynamics} (SINDy)
encompasses a variety of algorithms using different techniques to promote sparsity in $\ma C$,
many inspired by compressed sensing.
For example, Sequentially Thresholded Least Squares (STLS) (see, e.g.,~\cite{BD09})
can be used iteratively to remove small coefficients in $\ma C$~\cite{BPK16}.
Another approach is to add an $\ell_p$-norm regularization term to the objective, i.e.,
\begin{equation}\label{eq:lsoi_reg}
	\min_{\ma C} \| \ma D \ma Y - \ma \Phi(\ma Y) \ma C\|_\fro^2 + \lambda R_p(\vectorize (\ma C))
	\quad \text{where} 
	\quad R_p(\ve w) \coloneqq
	\begin{cases}
		\| \ve w\|_p^p, & p > 0; \\
		\| \ve w\|_0,  & p = 0;
	\end{cases}
\end{equation}
where $\vectorize(\cdot)$ denotes row-major vectorization
and $\|\ve w\|_0$ is the number of nonzero entries in $\ve w$.
This is a nonconvex problem for $p < 1$ and nondifferentiable for $p=0$;
however, in many applications small values of $p$, especially $p=0$,
provide better recovery.
There are a variety of techniques that can be used to solve~\cref{eq:lsoi_reg},
such as Iteratively Reweighted Least Squares (IRLS)~\cite[Subsec.~4.5.2]{Bjo96}
or Iteratively Reweighted $\ell_1$-norm (IR$\ell_1$)~\cite{CWB08}.
The latter is used for recovering dynamical systems in~\cite{CPD21}
with an accuracy better than that of STLS.
Our approach with SIDDS will be to promote sparsity
by adding a similar regularization penalty
\begin{equation}\label{eq:sidds_reg}
	\begin{split}
		\min_{\ma C, \ma Z} & \ \| \ma Y - \ma Z\|_\fro^2 + \alpha R_p(\vectorize(\ma C)) \\
		\text{s.t.} & \ \ma D \ma Z = \ma \Phi(\ma Z) \ma C.
	\end{split}
\end{equation}
We refer to this variant as SIDDS+$\ell_p$.
In our numerical experiments, we choose $p=0$
and use regularized IRLS~\cite{CY08} to provide a quadratic approximation of $R_p$.

Although these sparsity promoting techniques will often regularize the problem
and reduce the impact of noise,
simply identifying the correct sparsity structure is not sufficient
to obtain an accurate parameter estimate when using the LSOI objective.
Returning to the simple harmonic oscillator example of~\cref{eq:sho},
if we fix the correct sparsity structure for both of these methods, we obtain
\begin{center}
\vspace{0.5em}
\noindent
\begin{tabular}{@{}cc}
	LSOI + fixed sparsity & SIDDS + fixed sparsity \\[0.5em]
	{\small $ 
	\begin{bmatrix}
		\dot{x}_1 \\ \dot{x}_2 
	\end{bmatrix}
	= 
	\begin{bmatrix}
		0 &   7.15983089 \\
 		-5.69559158  & 0
	\end{bmatrix}
	\begin{bmatrix}
		x_1 \\ x_2
	\end{bmatrix}
	$}
	&
	{\small $ 
	\begin{bmatrix}
		\dot{x}_1 \\ \dot{x}_2
	\end{bmatrix}
	= 
	\begin{bmatrix}
		0 &  0.99695336 \\
		-1.01100718 &  0
	\end{bmatrix}
	\begin{bmatrix}
		x_1 \\ x_2
	\end{bmatrix}
	$}.
\end{tabular}
\vspace{0.5em}
\end{center}
Hence fixing the sparsity pattern---%
equivalent to hard thresholding in STLS---%
has not improved the parameter estimate using LSOI.
Thus sparsity promoting regularization is not sufficient to ensure
a method with an LSOI objective term will provide an accurate estimate.
Instead, noise must be removed from measurements $\ma Y$ 
so than an accurate estimate of $\ma C$ may be obtained.

\subsection{Existing Denoising Work}
There are a wide variety of algorithms for identifying sparse dynamics.
We focus our attention on those methods, like ours, 
that simultaneously estimate the coefficients and true state, or, equivalently, measurement noise.
For the most part, 
these methods do so by adding a penalty based on the dynamical system mismatch
rather than adding an equality constraint as with SIDDS.
Beyond these methods,
there are a variety of techniques that separately denoise measurements
and then apply LSOI or a related algorithm;
see, e.g.,~\cite{CPD22x}.
Although these yield improved performance in the presence of noise, 
their use requires careful parameter tuning.
Hence, we leave a more thorough comparison to future work.

\subsubsection{Sparse Corruption}
A related case to ours where where all measurements have been corrupted by noise,
is the case where an unknown (small) fraction of measurements are contaminated by noise.
If $\ma Y - \ma X$ has only a few nonzero rows, 
the dynamical system constraint will similarly be satisfied except for a few nonzero rows;
if $\ma D \ma X = \ma \Phi(\ma X) \ma C$, then the constraint mismatch $\ma E$ satisfies 
\begin{equation}
	\ma E= \ma \Phi(\ma Y) \ma C - \ma D\ma Y = \ma \Phi(\ma Y - \ma X) \ma C - \ma D (\ma Y - \ma X).
\end{equation}
Tran and Ward propose identifying these corrupted measurements
and a sparse dynamical system by solving~\cite[eq.~(9)]{TW17}
\begin{equation}
	\begin{split}
	\min_{\ma C\in \R^{n \times d},\  \ma E \in \R^{m \times d}} & \ \sum_{j} \| \ma E_{j,\cdot} \|_2 \\
	\text{s.t.} &\  \ma D \ma Y + \ma E = \ma \Phi(\ma Y) \ma C
		\quad \text{and $\ma C$ sparse}
	\end{split}
\end{equation}
where $\ma E_{j, \cdot}$ represents the $j$th row of $\ma E$.
This formulation aims for group-sparsity with respect to the rows of $\ma E$
by using a $\ell_1$ convex relaxation of the $\ell_0$-norm in the objective.
Although this method does identify the corrupted measurements,
it does not identify a corrected state
and thus avoids the nonlinearity in $\ma \Phi$, 
working with a fixed $\ma \Phi(\ma Y)$ rather than our $\ma \Phi(\ma Z)$.

\subsubsection{Modified SINDy}
Modified SINDy~\cite{KBK20x} 
introduces a variable $\ma N$ to estimate the noise
and penalizes violation of the dynamical system constraint by the denoised state.
Let $\op E^{t}(\ve x, \ma C)$ be the evolution operator 
that advances the differential equation with coefficients $\ma C$ from initial condition $\ve x$, 
\begin{equation}\label{eq:evolution}
	\op E^{t}: \R^d \times \R^{n \times d} \to \R^d, 
	\quad 
	\op E^{t}(\ve \zeta_0, \ma C) \coloneqq \ve \zeta(t),
	\quad \text{where} \quad
	\left\lbrace 
		\begin{aligned}
			\dve \zeta(t) &= \ve f(\ve \zeta(t); \ma C), \\
			\ve \zeta (0) &= \ve \zeta_0.
		\end{aligned}
	\right. 
\end{equation}
Modified SINDy then estimates $\ma C$ 
by adding a penalty to LSOI for mismatches in the denoised evolution $q$ steps forwards and backwards,
\begin{equation}
	\min_{\substack{ \ma C \in \R^{n \times d}\\ \ma N \in \R^{m \times d}}} \|\ma D (\ma Y - \ma N) - \ma \Phi(\ma Y - \ma N) \ma C \|_\fro^2 
		+ \sum_{j=q+1}^{M -q} \sum_{\substack{i=-q\\ i \ne 0}}^q
			\omega_i \| \ve y_{j+i} - \ve n_{j+i} - \op E^{\delta i }(\ve y_j - \ve n_j, \ma C)\|_2^2
\end{equation}
where $\omega_i>0$ is a weight.
Their implementation alternates between minimizing the objective
above using stochastic, gradient-based optimization 
and using STLS to identify a sparsity structure in $\ma C$;
a similar approach is used in~\cite{RKB19}.

Although more accurate than LSOI,
this approach is suboptimal because 
it imposes the dynamical system constraint through a penalty
rather than as an equality constraint.
We can observe the loss of accuracy  
in the simple harmonic oscillator example;
using their implementation of Modified SINDy,
we recover
\begin{center}
\vspace{0.5em}
\noindent
\begin{tabular}{@{}cc}
	Modified SINDy & SIDDS+$\ell_0$ \\[0.5em]
	{\small $ 
	\begin{bmatrix}
		\dot{x}_1 \\ \dot{x}_2 
	\end{bmatrix}
	=
	\begin{bmatrix}
		0 &  0.91451177 \\
		-1.10430015 &  0 
	\end{bmatrix}
	\begin{bmatrix}
		x_1 \\ x_2
	\end{bmatrix}
	$}
	&
	{\small $ 
	\begin{bmatrix}
		\dot{x}_1 \\ \dot{x}_2 
	\end{bmatrix}
	= 
	\begin{bmatrix}
		0 & 0.99695336 \\
 		-1.01100718 &  0       
	\end{bmatrix}
	\begin{bmatrix}
		x_1 \\ x_2
	\end{bmatrix}
	$}.
\end{tabular}
\vspace{0.5em}
\end{center}
Although Modified SINDy vastly outperforms LSOI, 
the coefficient error is an order of magnitude larger than SIDDS.

\subsubsection{Physics Informed Spline Learning}
Physics Informed Spline Learning (PiSL)~\cite{SLS21}
uses a similar approach to Modified SINDy,
adding a penalty for violating the dynamical system constraint.
In PiSL, 
the estimated state $\ve z(t)$ is expressed in a cubic spline basis,
$\ve z(t) = \ma T(t) \ma P$
where $\ma T(t)$ is a cubic spline basis and $\ma P$ are the control points.
To identify the dynamical system, PiSL solves~\cite[eq.~(9)]{SLS21}
\begin{equation}
	\min_{\ma P, \ma C}
		\frac{1}{m} \sum_{j=1}^m \big\| \ve y_j - \ma T(t_j)\ma P\big\|_2^2
		+ \frac{\alpha}{m}
			\sum_{j=1}^m \big\| \dma T(t_j) \ma P - \sum_{k=1}^n \ve c_k \phi(\ma T(t_j) \ma P)\big\|_2^2 
		+\beta\| \ma C\|_0
\end{equation}
(we omit the stochastic subsampling used in the first two terms in the objective).
PiSL's implementation uses a similar approach to Modified SINDy:
a minimization alternating between $\ma P$ and $\ma C$ with a fixed sparsity structure
and then applying STLS to identify the sparsity pattern in $\ma C$.
Although this method uses a different basis
and differs in some details, 
the overall approach is similar to Modified SINDy.


\subsection{Overview}
The remainder of this manuscript is answers to two questions:
what is the statistical performance of SIDDS 
and how can we efficiently solve the SIDDS optimization problem?
First,
we convert matrix quantities in~\cref{eq:sidds} to vector quantities
simplify analysis in \cref{sec:notation}.
Then in \cref{sec:crb} we obtain a lower bound on the covariance
of any dynamical system estimator using the constrained Cram\'er-Rao Lower Bound (CRLB).
We also derive asymptotic estimates of bias and covariance of both SIDDS and LSOI;
numerical experiments show SIDDS obtains the CRLB up to discretization error
whereas LSOI does not.
Next in \cref{sec:algorithm} we show how to efficiently implement SIDDS
using an IRLS approximation to the $\ell_p$-norm regularization 
and a sequential quadratic program (SQP) with a preconditioned MINRES iteration.
Finally in \cref{sec:examples},
we provide a number of numerical experiments
comparing the performance of SIDDS and SIDDS+$\ell_0$
to LSOI, SINDy with STLS, and Modified SINDy.
These experiments show that SIDDS
almost exactly obtains the CRLB 
in a variety of settings
whereas existing techniques do not.
Moreover, SIDDS+$\ell_0$ accurately recovers the sparsity structure of $\ma C$ 
at higher levels of noise than other methods.
We conclude with a brief discussion
of future directions in \cref{sec:discussion}.

\subsection{Reproducibility}
Following the principles of reproducible research,
software implementing our algorithms 
and the scripts to generate data in the figures
are available at {\tt \url{http://github.com/jeffrey-hokanson/sidds}}.

\section{Notation and Derivatives\label{sec:notation}}
To follow standard practice in optimization, 
we reformulate the SIDDS optimization problem
with a vector-valued objective and constraint
which simplifies the derivations that follow.
In this section, 
we introduce notation for these vectorized quantities,
the constraint and its derivative,
and generalize the SIDDS optimization problem.

\subsection{Vectorization}
Throughout, we use the row-major vectorization operator
\begin{equation}
	\vectorize: \R^{m\times d} \to \R^{md}, \qquad 
	\vectorize(\ma X) = 
	\vectorize\left(
		\begin{bmatrix} \ve x_1^\trans \\ \vdots \\ \ve x_m^\trans \end{bmatrix}
	\right) 
	\coloneqq \begin{bmatrix} \ve x_1 \\ \vdots \\ \ve x_m\end{bmatrix}.
\end{equation}
For brevity, vectorized matrices
are denoted by the corresponding lower case letter annotated with a harpoon; e.g.,
\begin{align}
	\vve x &\coloneqq \vectorize(\ma X), &
	\vve y &\coloneqq \vectorize(\ma Y), &
	\vve z &\coloneqq \vectorize(\ma Z), &
	\vve c &\coloneqq \vectorize(\ma C).
\end{align}
Quantities that are still matrices after vectorization
are denoted by uppercase letters annotated with a harpoon; e.g., 
\begin{align}
	\vectorize(\ma \Phi(\ma Z)\ma C) &= 
	\vma \Phi(\vve z) \vve c, &
	\vma \Phi(\vve z) &\coloneqq \ma \Phi(\vve z) \otimes \ma I_d, \\
	\vectorize(\ma D \ma Z) &=
		\vma D \vve z, &
	\vma D &\coloneqq \ma D \otimes \ma I_d, 
\end{align}
where $\otimes$ denotes the Kronecker product
and $\ma I_d\in \R^{d\times d}$ denotes the identity.
The choice of row-major vectorization is important
so that $\vma D$ has low bandwidth.
If $\ma D$ uses a $q$-point stencil, 
then $\vma D$ has bandwidth $d(q-1)$.

\subsection{Constraint}
The constraint and its derivative are:
\begin{align}\label{eq:constraint}
	\vve h(\vve c, \vve z) \coloneqq& 
		\vectorize(\ma D \ma Z - \ma \Phi(\ma Z)\ma C)
		= \vma D \vve z - \vma \Phi(\vve z)\vve c \\
	\label{eq:constraint_der}
	\nabla \vve h(\vve c, \vve z) = & 
	\begin{bmatrix} 
		\nabla_{\vve c} \vve h(\vve c, \vve z) & \nabla_{\vve z} \vve h(\vve c, \vve z)
	\end{bmatrix}
	= 
	\begin{bmatrix}
		- \vma \Phi(\vve z) & \vma D - \nabla_{\vve z}[\vma \Phi(\vve z)\vve c]
	\end{bmatrix}.
\end{align}
The last term merits some elaboration.
The gradient of $\vma \Phi$ is
the 3-tensor $\nabla \vma \Phi: \R^{dm} \to \R^{(dm)\times (dn) \times (dm)}$;
thus a Taylor series expansion of $\vma \Phi$ about $\vve z$ is
\begin{equation}
	\vma \Phi(\vve z + \vve {\delta z})
	= \vma \Phi(\vve z)
	+ \nabla \vma \Phi(\vve z) \tenmult3 \vve {\delta z}
	+ \order( \| \vve {\delta z}\|_2^2),
\end{equation}
where $\tenmult3$ denotes tensor multiplication along 
the third tensor mode~\cite[subsec.~2.5]{KB09}.
If we take the same Taylor expansion for $\vma \Phi(\vve z) \vve c$,
\begin{align}
	\vma \Phi(\vve z + \vve {\delta z})\vve c
	&= \vma \Phi(\vve z)\vve c
	+ \left[\nabla \vma \Phi(\vve z) \tenmult3 \vve {\delta z}\right]\vve c
	+ \order( \| \vve {\delta z}\|_2^2), \\
	\intertext{we can interchange the order of tensor multiplication yielding}
	\vma \Phi(\vve z + \vve {\delta z})\vve c
	&= \vma \Phi(\vve z)\vve c
	+ \left[\nabla \vma \Phi(\vve z) \tenmult2 \vve c \right]\vve {\delta z}
	+ \order( \| \vve {\delta z}\|_2^2).
\end{align}
Hence the derivative of $\vma \Phi(\vve z) \vve c$ with respect to $\vve z$ is
\begin{equation}
	\nabla_{\vve z} [\vma \Phi(\vve z)\vve c] = \nabla \vma \Phi(\vve z) \tenmult2 \vve c \in \R^{(dm)\times(dm)}.
\end{equation}
This expression allows us to reduce the memory required to compute the constraint derivative;
our implementation never explicitly forms the large 3-tensor $\nabla \vma \Phi(\vve z)$,
but instead builds $\nabla \vma \Phi(\vve z) \tenmult2 \vve c$ directly.

\subsection{Problem Statement\label{sec:notation:problem}}
Using this vectorized notation, we can restate and generalize 
the SIDDS+$\ell_p$ problem~\cref{eq:sidds_reg} as
\begin{equation}\label{eq:sidds_weight}
	\begin{split}
		\min_{\vve c \in \R^{dn}, \vve z\in \R^{dm}} \ & 
			(\vve y- \vve z)^\trans \vma M (\vve y - \vve z)  + \alpha R_p(\vve c)  \\
		\text{s.t.} \ & \vve h(\vve c, \vve z) = \vma D\vve z - \vma \Phi(\vve z) \vve c = \ve 0.
	\end{split}
\end{equation}
Here the symmetric positive definite matrix $\vma M$ 
acts as a weight on the mismatch between the measurements and estimated state.
There are two important situations when we might use $\vma M \ne \ma I$.
If $\vve y \sim \mathcal{N}(\vve x, \vma \Sigma)$,
then taking $\vma M = \vma \Sigma^{-1}$ yields a near optimal estimate for $\vve c$ and $\vve z$
as shown in \cref{sec:crb:sidds}.
The other situation is when we need to decouple the time step used in the ODE constraint
from sample rate of measurements in order to improve recovery as illustrated in \cref{sec:examples:integration}.
\section{Statistical Performance\label{sec:crb}}
A key question is:
how do methods like SIDDS and LSOI
perform when measurements are contaminated by noise?
In this section, we analyze these methods when measurements $\vve y$
have been contaminated with normally distributed additive noise $\vve n$ with zero mean and
full-rank covariance $\vma \Sigma$,
\begin{equation}\label{eq:noise_model}
\vve y = \vve x + \vve n, \qquad \vve n \sim \mathcal{N}(\vve 0, \vma \Sigma);
\quad \text{equivalently,} \quad
\vve y \sim \mathcal{N}(\vve x, \vma \Sigma).
\end{equation}
For these methods, we estimate how $\vve c$ is perturbed in the limit of small noise
and compute the mean and covariance of $\vve c$ asymptotically.
Although only asymptotically valid, 
numerical experiments 
show these estimates provide reliable guides to performance.
To begin, we compute a lower bound on the covariance of $\vve c$
based on the solution to the continuous optimization problem.
We then compare the asymptotic covariance estimates for LSOI and SIDDS to this lower bound;
LSOI is far larger than this lower bound, whereas SIDDS satisfies it almost exactly.

\subsection{Covariance Lower Bound}
Here we use the constrained variant of the Cram\'er-Rao lower bound 
to bound the covariance of $\vve c$.
To begin, 
let us consider the likelihood function associated with the noise model in~\cref{eq:noise_model}.
Suppose our model has estimated the system state sequence $\vve z$ using coefficients $\vve c$;
the corresponding likelihood is
\begin{equation}\label{eq:likelihood}
	p(\vve y; \vve c, \vve z) \coloneqq
	\det(2\pi \vma \Sigma)^{-\frac12} 
	\exp[ -\tfrac12 (\vve y - \vve z)^\trans \vma \Sigma^{-1} (\vve y- \vve z)];
\end{equation}
see, e.g.,~\cite[Sec.~2.2]{SW89}.
This likelihood function does not depend on $\vve c$ explicitly;
however $\vve z$ should satisfy the evolution equations,
$\ve z_{j+1} = \op E^{\delta j}(\ve z_1, \ma C)$ (recall $\op E$ was defined in~\cref{eq:evolution}).
We denote this constraint set as
\begin{equation}
	\Omega
		\coloneqq \lbrace \vve c, \vve z: \op E^{\delta j}(\ve z_1, \vve c) =  \ve z_{j+1}, \ j=1,\ldots, m-1
		\rbrace.
\end{equation}
We compute the constrained Cram\'er-Rao lower bound following~\cite{BE09}
using the Fisher information matrix $\vma J$ 
and an orthogonal basis for the tangent space of constraints $\vma U(\vve c, \vve z)$.
For this likelihood function~\cref{eq:likelihood}, the Fisher information matrix is constant,
\begin{equation}
	\vma J \coloneqq \E_{\vve y} \left[
		\nabla^2_{\vve c, \vve z} \log p(\vve y; \vve c, \vve z) 
	\right]
	= 
	\begin{bmatrix}
		\ma 0 & \ma 0 \\ \ma 0 & \vma \Sigma^{-1}
	\end{bmatrix}
\end{equation}
where $\E_{\vve y}$ denotes the expectation over $\vve y$.
Rewriting the constraint set $\Omega$ in terms of a function 
$\vve h_\Omega: \R^{dn} \times \R^{dm} \to \R^{d(m-1)}$
\begin{equation}
	\vve h_\Omega(\vve c, \vve z) \coloneqq 
		\begin{bmatrix} 
			\op E^{\delta}( \ve z_1, \vve c) - \ve z_2 \\ 
			\vdots \\ 
			\op E^{\delta (m-1)} (\ve z_1, \vve c) - \ve z_m
		\end{bmatrix},
\end{equation}
the tangent space of the constraint 
$\vma U(\vve c, \vve z) \in \R^{(d(m + n)) \times (d(n+1))}$ satisfies
\begin{equation}
	\ma 0 = [\nabla_{\vve c, \vve z} \vve h_\Omega(\vve c, \vve z) ] \vma U(\vve c, \vve z),
	\quad \text{and} \quad
	\vma U(\vve c, \vve z)^\trans \vma U(\vve c, \vve z) = \ma I.
\end{equation}
Let $\vve c^\star$ and $\vve z^\star = \vve x$ be the true coefficients and states respectively.
Then for any estimator that given $\vve y$ produces unbiased estimates
$\vve c^\generic(\vve y)$ 
and $\vve z^\generic(\vve y)$,
namely, 
\begin{equation}
	\E_{\vve y} \vve c^\generic(\vve y)=\vve c^\star \quad  \text{and}  \quad \E_{\vve y} \vve z^\generic(\vve y) =\vve z^\star = \vve x,
\end{equation}
then $\vve c^\generic$ and $\vve z^\generic$ satisfy the constrained Cram\'er-Rao lower bound
\begin{equation}\label{eq:crlb}
	\Cov_{\vve y} \begin{bmatrix}\vve c^\generic(\vve y) \\ \vve z^\generic(\vve y) \end{bmatrix}
	\succeq 
	\vma U [\vma U^\trans \vma J \vma U]^{+} \vma U^\trans,
	\quad \vma U = \vma U(\vve c^\star, \vve z^\star),
\end{equation}
where $\Cov_{\vve y}$ denotes the covariance with respect to $\vve y$, 
$^+$ the pseudoinverse~\cite[Sec.~5.5.2]{GL13},
and $\succeq$ the ordering of positive semidefinite matrices~\cite[sec.~7.7]{HJ12}. 

\subsubsection[Nuisance Variables]{Covariance of $\vve c$}
As our goal is to estimate the coefficients $\vve c$, 
the estimated state $\vve z$ is considered a nuisance variable.
To compute the covariance of $\vve c$ alone,
we use the selection matrix $\vma S_{\vve c}\in \R^{(d(m+n)) \times n_c}$ 
to pick those columns corresponding to $\vve c$ where $n_c$ is the number of entries in $\vve c$
($n_c \ne dn$ with a sparsity constraint on $\vve c$).
Here, $\ma S_{\vve c}$ picks the $(1,1)$ block of the CRLB: 
\begin{equation}\label{eq:crlb_c}
	\Cov \vve c^\generic \succeq 
		\vma S_{\vve c}^\trans \vma U [\vma U^\trans \vma J \vma U]^+ \vma U^\trans \vma S_{\vve c}.
\end{equation}

\subsubsection{Sparsity Constraint\label{sec:crb:sparse}}
Imposing a sparsity structure on $\vve c$ changes the Cram\'er-Rao lower bound
as this adds an additional constraint.
Suppose $\set I(\vve c)$ selects indices of $\vve c$;
then setting these entries to zero is equivalent to the constraint set
\begin{equation}
	\Omega_{\text{sparse}} = \Omega \cap \lbrace \vve c, \vve z: \set I(\vve c) = \ve 0 \rbrace.
\end{equation}
As this is a larger set of constraints,
this decreases the size of the tangent space.
In general, 
this yields a smaller covariance of the remaining nonzero entires.


\subsubsection{Computing the Tangent Space}
We can compute an orthogonal basis $\vma U$ for the tangent space of $\Omega$
by solving the sensitivity equations.
Let $\ma V$ be the Jacobian with respect to the coefficients $\vve c$
and $\ma W$ be the Jacobian with respect to the initial conditions
\begin{align}
	\ma V(t) &= \nabla_{\vve c} \op{E}^t(\ve z_1, \vve c) \in \R^{d \times n_c},
	&
	\ma W(t) &= \nabla_{\ve z_1} \op{E}^t(\ve z_1, \vve c) \in \R^{d\times d}.
\end{align}
The vectorized versions of these two quantities, $\vve v(t)$ and $\vve w(t)$,
evolve according to the coupled differential equations 
\begin{equation}
	\left\lbrace
	\begin{alignedat}{2}
		\partial_t  \ve z(t) &= \ve f(\ve z(t), \vve c) 
	& \ve z(0) &= \ve z_0\\
		\partial_t \vve v(t) &= \ma F(\ve z(t), \vve c) \vve v(t)
			+ \nabla_{\vve c} \ve f(\ve z(t), \vve c)
		\qquad & \vve v(0) &= \ve 0 \\
		\partial_t \vve w(t) &=
			\ma F(\ve z(t), \vve c) \vve w(t) 
		& \vve w(0) &= \vectorize(\ma I_d) 
	\end{alignedat}
	\right.
\end{equation}
where $\ma F(\ve z,\vve c) = \nabla_{\ve z} \ve f(\ve z, \vve c) = \sum_{k} \ve c_k \nabla_{\ve z} \phi_k(\ve z)$.
Thus the gradient of constraints $\vve h_\Omega$ is
\begin{align}
	\setlength\arraycolsep{2pt} 
	\nabla_{\vve c, \vve z}
	\vve h_\Omega(\vve c, \vve z)
	= 
	\nabla_{\vve c, \vve z}
	\begin{bmatrix}
		\op E^{\delta}(\ve z_1, \vve c) - \ve z_2\\
		\op E^{2\delta}(\ve z_1, \vve c) - \ve z_3\\
		\vdots \\
		\op E^{m\delta}(\ve z_1, \vve c) - \ve z_m
	\end{bmatrix}
	\! = \!
	\begin{bmatrix}
		\ma V(\delta) & \ma W(\delta) & - \ma I	\\
		\ma V(2\delta) & \ma W(2\delta) & & -\ma I \\
		\vdots & \vdots & & & \ddots \\
		\ma V(m\delta) & \ma W(m \delta) & & & & -\ma I
	\end{bmatrix}\!.
\end{align}
The structure in this matrix allows
us to write down an explicit formula its nullspace,
\begin{equation}
	\Nullspace \nabla_{\vve c, \vve z} \vve h_\Omega(\vve c, \vve z)= 
	\Range 
	\left(
	\begin{bmatrix}
		\ma I & \ma 0 \\
		\ma 0 & \ma I \\
		\ma V(\delta) & \ma W(\delta) \\
		\vdots & \vdots \\
		\ma V(m\delta) & \ma W(m \delta) 
	\end{bmatrix}
	\right).
\end{equation}
Thus we can compute an orthogonal basis for the tangent space of $\Omega$ at $\vve c, \vve z$,
namely, $\vma U(\vve c, \vve z)$,
by performing a reduced QR-factorization of this matrix above.

\subsection{Performance of LSOI\label{sec:crb:lsoi}}
Let $\vve c^\lsoi$ be the coefficient estimate of LSOI~\cref{eq:lsopinf}, 
\begin{equation}
	\vve c^\lsoi \coloneqq \argmin_{\vve c} \|\vma D \vve y - \vma \Phi(\vve y )\vve c\|_2^2.
\end{equation}
We can alternatively write this solution using the pseudoinverse $^+$,
\begin{equation}
	\vve c^\lsoi = \vve c^\lsoi(\vve y) = \vma \Phi^+(\vve y) \vma D\vve y.
\end{equation}
Assuming $\vve y= \vve x + \sigma \vve n$ with $\vve n \sim \mathcal{N}(\ve 0 , \ma \Sigma)$,
in the limit of small noise ($\sigma \to 0$), 
\begin{align}
	\vve c^\lsoi &= \vma \Phi^+(\vve x + \sigma\vve n)\vma D (\vve x + \sigma\vve n) \\
		&= \vma \Phi^+(\vve x) \vma D \vve x 
		+ \sigma \left([\nabla \vma \Phi^+(\vve x)]\tenmult3 \vve n \right)\vma D \vve x
		+ \sigma \vma \Phi^+(\vve x)\vma D \vve n + \order(\sigma^2) \\
		&= \vma \Phi^+(\vve x) \vma D \vve x 
		+ \sigma \left([\nabla \vma \Phi^+(\vve x)]\tenmult2 (\vma D \vve x) + \vma \Phi^+(\vve x) \vma D\right)\vve n
		+ \order(\sigma^2).
\end{align}
Let this first order transformation of $\vve n$ be denoted by
\begin{equation}
		\vma T_\lsoi \coloneqq 
		[\nabla \vma \Phi^+(\vve x)]\tenmult2 (\vma D \vve x) + \vma \Phi^+(\vve x) \vma D.
\end{equation}
To first order, $\vve c^{\lsoi}$ has a bias proportional to finite difference error in the derivative,
\begin{align}
	\E [\vve c^\lsoi - \vve c^\star] &= \vma \Phi^+(\vve x) \vma D \vve x - \vve c^\star + \E_{\vve n} [\vma T_\lsoi \vve n] + \order(\sigma^2) \\
		 &= \vma \Phi^+(\vve x) (\vma D \vve x - \dot{\vve x})   + \order(\sigma^2),
\end{align}
since $\dot{\vve x} = \vma \Phi(\vve x)\vve c^\star$
where $\dot{\vve x}$ denotes the state derivative.
Similarly, the covariance is
\begin{equation}\label{eq:cov_lsoi}
	\Cov \vve c^\lsoi = \sigma^2 \vma T_\lsoi \vma \Sigma \vma T_\lsoi^\trans + \order(\sigma^3).
\end{equation}
We can explicitly compute $\vma T_\lsoi$ by using Golub and Pereyra's formula
for the derivative of a pseudoinverse~\cite[eq.~(4.12)]{GP73};
denoting by $\partial_i$ the derivative of the $i$th entry of~$\vve c$
and omitting arguments,
\begin{align}
	\partial_i \vma \Phi^+\!
	= \! - \vma \Phi^+ [\partial_i \vma \Phi] \vma \Phi^+
	\! + \! \vma \Phi^+ \vma \Phi^{+\trans} [\partial_i \vma \Phi^\trans](\ma I \! - \! \vma \Phi \vma \Phi^+)
	\! + \! (\ma I\! - \! \vma \Phi^+ \vma \Phi)[\partial_i \vma \Phi^\trans] \vma \Phi^{+\! \trans} \vma \Phi^+.
\end{align}

Next, using this asymptotic analysis,
we demonstrate LSOI yields suboptimal estimates.
In the following examples, 
we consider the simple harmonic oscillator example from the introduction
with a fixed sparsity pattern selecting the true, nonzero values of $\vve c$.
Taking $\vve n\sim \mathcal{N}(\vve 0, \sigma^2\ma I)$ with $\sigma=10^{-2}$, 
we compute a Monte Carlo estimate of the covariance of $\vve c^\lsoi$,
the asymptotic estimate~\cref{eq:cov_lsoi},
and the Cram\'er-Rao lower bound~\cref{eq:crlb_c}
\begin{center}
\vspace{0.5em}
\noindent
\begin{tabular}{@{}ccc}
	Monte Carlo $\Cov[\vve c^\lsoi]$ & Asymptotic $\Cov[\vve c^\lsoi]$ & Unbiased CRLB  \\[0.5em]
	{\small$\sigma^2 \writemat{data/fig_sho_cov_lsopinf.csv}{} $}
	&
	{\small$\sigma^2 \writemat{data/fig_sho_cov_analytic_lsoi_pts3.csv}{} $}
	&
	{\small $\sigma^2\writemat{data/fig_sho_cov_analytic_crb.csv}{}$}.
\end{tabular}
\vspace{0.5em}
\end{center}
This example shows 
that the asymptotic estimate provides an accurate covariance
and the covariance LSOI is substantially larger than the lower bound.
This example used a 3-point finite difference rule;
we might expect ahigher order finite difference approximation would decrease the covariance.
In fact, the opposite happens!
Examining the asymptotic covariance for this same problem, we observe
\begin{center}
\vspace{0.5em}
\noindent
\begin{tabular}{@{}ccc}
	3-point rule $\Cov[\vve c^\lsoi]$ & 5-point rule $\Cov[\vve c^\lsoi]$ & 7-point rule $\Cov[\vve c^\lsoi]$ \\[0.5em]
	{\small $\sigma^2\writemat{data/fig_sho_cov_analytic_lsoi_pts3.csv}{}$ }
	&
	{\small $\sigma^2\writemat{data/fig_sho_cov_analytic_lsoi_pts5.csv}{}$}
	&
	{\small $\sigma^2\writemat{data/fig_sho_cov_analytic_lsoi_pts7.csv}{}$.}
\end{tabular}
\vspace{0.5em}
\end{center}
The exact origin of this effect is unclear,
but it is likely a result of the increased bandwidth of $\vma D$.
This result provides yet another reason that high order differencing schemes 
are not frequently seen in other work.

\subsection{Performance of SIDDS\label{sec:crb:sidds}}
Here we consider the SIDDS problem with a weighted objective
as introduced in~\cref{eq:sidds_weight}.
We will show that by choosing $\vma M = \vma \Sigma^{-1}$,
we obtain estimates that approximately satisfy the CRLB.

As before, we consider the limit of small noise $\vve y = \vve x + \sigma \vve n$
with $\sigma \to 0$ and construct perturbation estimates around the true values
\begin{align}
	\vve c^\sidds &= \vve c^\star + \sigma \vve c^{(1)} + \order(\sigma^2)& 
	\vve z^\sidds &= \vve x + \sigma \vve z^{(1)} + \order(\sigma^2). 
\end{align}
Linearizing the constraint of~\cref{eq:sidds_weight}
around the zeroth order terms yields the quadratic program
with Karush-Kuhn-Tucker (KKT) system
\begin{align}
	\begin{bmatrix}
		\ve 0 & & \vma K^\trans \\
		&  \sigma^2 \vma M & \vma L^\trans \\
		\vma K & \vma L & \ve 0
	\end{bmatrix}
	\begin{bmatrix}
		\vve c^{(1)} \\ \vve z^{(1)} \\ \vve w
	\end{bmatrix}
	=
	\begin{bmatrix}
		\vve 0 \\ \sigma^2 \vma M \vve n \\ -\vve h(\vve c, \vve x)
 	\end{bmatrix}
\end{align}
where $\vve h$ is the constraint defined in~\cref{eq:constraint}
and $\vma K$ and $\vma L$ are the two blocks in the constraint derivative~\cref{eq:constraint_der}.
To compute $\vve c^{(1)}$ and $\vve h^{(1)}$ we use a reduced Hessian approach 
where $\vma U_\sidds$ is an orthogonal basis for the nullspace of $\begin{bmatrix} \vma K & \vma L\end{bmatrix}$,
\begin{align}\label{eq:sidds:first_order}
	\begin{bmatrix} \vve c^{(1)} \\ \vve z^{(1)} \end{bmatrix}
	= -\begin{bmatrix} \vma K & \vma L\end{bmatrix}^+ \vve h(\vve c^\star, \vve x) 
		+	 \vma U_\sidds \left(
			\vma U_\sidds^\trans
			\begin{bmatrix} 
				\ve 0 & \\
				&  \vma M
			\end{bmatrix}
			\vma U_\sidds
		\right)^{+} \vma U^\trans_\sidds \begin{bmatrix} \vve 0 \\ \vma M \vve n \end{bmatrix}.
\end{align}
As with LSOI, we note that SIDDS has a slight bias since 
the discretized constraint $\vve h(\vve c^\star, \vve x)$ is not necessarily satisfied exactly
\begin{equation}\label{eq:sidds:bias}
	\E \begin{bmatrix} 
		\vve c^\sidds - \vve c^\star \\
		\vve z^\sidds - \vve x
	\end{bmatrix}
	= -\begin{bmatrix} \vma K & \vma L\end{bmatrix}^+ \vve h(\vve c^\star, \vve x) + \order(\sigma^2).
\end{equation}
Using a higher order finite difference scheme reduces this bias.
Returning to the simple harmonic oscillator example with fixed sparsity
and $\vve n \sim \mathcal{N}(\vve 0, \vve I)$,
\begin{center}
\vspace{0.5em}
\noindent
\begin{tabular}{@{}ccc}
	3-point $\E[\vve c^\sidds - \vve c^\star]$ & 
	5-point $\E[\vve c^\sidds - \vve c^\star]$ & 
	7-point $\E[\vve c^\sidds - \vve c^\star ]$ \\[0.5em]
	{\small$\writemat{data/fig_sho_cov_analytic_lide_pts3_bias.csv}{} $}
	&
	{\small$\writemat{data/fig_sho_cov_analytic_lide_pts5_bias.csv}{} $}
	&
	{\small$\writemat{data/fig_sho_cov_analytic_lide_pts7_bias.csv}{} $}.
\end{tabular}
\vspace{0.5em}
\end{center}
Note that this bias almost exactly matches the error in the example in~\cref{sec:intro:sparse}.

We obtain the asymptotic covariance of the SIDDS estimates 
by taking the expectation of the outer product of the first order perturbation~\cref{eq:sidds:first_order}, 
\begin{equation}
	\Cov \begin{bmatrix} \vve c^\sidds \\ \vve z^\sidds \end{bmatrix}
	= 
		\vma U_\sidds \! \left(
			\vma U_\sidds^\trans
			\begin{bsmallmatrix} 
				\ve 0 & \\
				&  \vma M
			\end{bsmallmatrix}
			\vma U_\sidds
		\right)^{\!\!+} \vma U^\trans_\sidds 
		\begin{bsmallmatrix} \vma 0 & \\ & \sigma^2 \vma M \vma \Sigma \vma M^\trans \end{bsmallmatrix}
		\vma U_\sidds \! \left(
			\vma U_\sidds^\trans
			\begin{bsmallmatrix} 
				\ve 0 & \\
				&  \vma M
			\end{bsmallmatrix}
			\vma U_\sidds
		\right)^{\!\!+} \vma U^\trans_\sidds 
	+\order(\sigma^3).
\end{equation}
If we take $\vma M = \vma \Sigma^{-1}$,
this subtantially simplifies:
\begin{equation}
	\Cov \begin{bmatrix} \vve c^\sidds \\ \vve z^\sidds \end{bmatrix}
	= 
		\vma U_\sidds \left(
			\vma U_\sidds^\trans
			\begin{bmatrix} 
				\ve 0 & \\
				&  \vma \Sigma^{-1}
			\end{bmatrix}
			\vma U_\sidds
		\right)^{+}
		\vma U_\sidds^\trans
	+\order(\sigma^3),
	\quad \text{if} \quad \vma M = \vma \Sigma^{-1}.
\end{equation}
This has the same form as \cref{eq:crlb},
except the basis for the true tangent space of the ODE constraint $\vma U$
has been replaced with that of the discretized ODE constraint $\vma U_\sidds$.
Typically, the subspace angle between these two spaces is small
leading the SIDDS estimate $\vve c^\sidds$ to nearly obtain the lower bound.
In the simple harmonic oscillator example with fixed sparsity,
\begin{center}
\vspace{0.5em}
\noindent
\begin{tabular}{@{}ccc}
	3-point $\Cov \vve c^\sidds$ & 5-point $\Cov \vve c^\sidds$ & Unbiased CRLB \\[0.5em]
	{\small$\writemat{data/fig_sho_cov_analytic_lide_pts3.csv}{} $}
	&
	{\small$\writemat{data/fig_sho_cov_analytic_lide_pts5.csv}{} $}
	&
	{\small $\writemat{data/fig_sho_cov_analytic_crb.csv}{}$}.
\end{tabular}
\vspace{0.5em}
\end{center}
Unlike LSOI, higher order difference schemes for SIDDS yield better estimates.
Note that the 3-point covariance is slightly smaller than the unbiased Cram\'er-Rao lower bound
which is possible due to its bias;
the 5-point covariance matches the lower bound with an error around $10^{-7}$.

\section{Solving SIDDS\label{sec:algorithm}}
We now turn our attention to devising an efficient numerical algorithm
to solve the SIDDS optimization problem:
\begin{equation}\label{eq:sidds_opt}
	\begin{split}
		\min_{\vve c,\,  \vve z } \ &  
			(\vve y - \vve z)^\trans \vma M (\vve y - \vve z) + \alpha R_p(\vve c)\\ 
		\text{s.t.} \ & \vma D \vve z = \vma \Phi(\vve z) \vve c;
	\end{split}
	\quad \text{recall} \quad 
	R_p(\vve c) \coloneqq
	\begin{cases}
		\| \vve c\|_p^p, & p > 0; \\
		\| \vve c\|_0, & p =0;
	\end{cases}
\end{equation}
and $\vma M \in \R^{(dm)\times(dm)}$ is positive semidefinite.
Structurally, this problem is close to a quadratic program,
and would be if $\vma \Phi(\vve z)$ was held constant and regularization removed ($\alpha=0$).
Hence, we will use a SQP approach to solve SIDDS.
However, rather than optimizing~\cref{eq:sidds_opt} directly, 
we will use IRLS to provide a convex approximation of $R_p$
as described in \cref{sec:algorithm:irls}.
We then use the SQP approach of Lin and Yuan~\cite{LY11}
to optimize the IRLS approximated problem,
taking special care to exploit structure;
we discuss this in \cref{sec:algorithm:sqp}.
We then show how this approach can be extended to problems
where multiple trajectories in \cref{sec:optimization:traj}.
We begin by discussing how initialize this nonconvex optimization problem.

\subsection{Initialization}
Since~\cref{eq:sidds_opt} is a nonconvex optimization problem, 
convergence to a global optimum from arbitrary initial estimates is not guaranteed.
With small noise, simply taking $\vve z^{(0)} = \vve y$ and $\vve c^{(0)}=\ve 0$
(superscripts denoting iteration number)
yields good local minimizers.
With large noise, 
we find better solutions by applying modest amount of smoothing.
In particular, we use Tikhonov smoothing,
\begin{equation}\label{eq:smooth}
	\vve z^{(0)} = \argmin_{\vve z} 
	\left\| 
		\begin{bmatrix}  \vve y \\ \ve 0 \end{bmatrix}
		-
		\begin{bmatrix} \ma I \\ \lambda \vma D^2\end{bmatrix} \vve z
	\right\|_2^2
\end{equation} 
where $\lambda>0$ is a smoothing parameter (in our experiments $\lambda=10^{-2}$)
and $\vma D^2$ is a 3-point finite-difference approximation of the second derivative.
Then we apply LSOI to estimate $\vve c^{(0)}$ based on $\vve z^{(0)}$.

\subsection{Approximating Regularization via IRLS\label{sec:algorithm:irls}}
The regularization $R_p$ in the objective~\cref{eq:sidds_opt}
presents a challenge:
for $p<1$, $R_p$ is concave
introducing negative curvature into the Hessian of the objective.
To avoid this difficulty, 
we replace $R_p$ with an IRLS approximation~\cite[subsec.~4.5.2]{Bjo96}
yielding a positive semidefinite Hessian.
At the $\ell$th iterate, the IRLS approximation of $R_p$ is
\begin{equation}
	R_p(\vve c^{(\ell)})
		\coloneqq \sum_{i} |c_i^{(\ell)} |^{p}
		= \sum_{i} |c_i^{(\ell)} |^{p-2}|c_i^{(\ell)}|^{2}
		\approx \sum_i
		|c_i^{(\ell - 1)}|^{p-2} |c_i^{(\ell)}|^2.
\end{equation}
If $\vve c^{(\ell-1)}\to \vve c^{(\ell)}$ as $\ell\to \infty$,
this IRLS approximation approaches $R_p$.
Numerically, when $c_i^{(\ell-1)}$ is small and $p<2$, we can encounter a divide by zero error.
Following Chartrand and Yin~\cite{CY08},
we avoid this by introducing an iteration dependent regularization $\epsilon^{(\ell)}>0$;
this yields the quadratic approximation:
\begin{equation}
	R_p(\vve c) \approx
	R_p^{(\ell)}(\vve c)
		\coloneqq
		\vve c^\trans \vma W^{(\ell-1)} \vve c,
	\qquad
	\vma W^{(\ell)} \coloneqq \diag
		\left(
		\left\lbrace
		[ \;  |c_i^{(\ell)}|^2 + \epsilon^{(\ell)}]^{\frac{p}{2} -1}
		\right\rbrace_i
	\right).
\end{equation}  
As the optimization proceeds, we let $\epsilon^{(\ell)}\to 0$.
We opt for the simple heuristic of decreasing $\epsilon$ 
by a factor of ten when optimization terminates successfully following~\cite{CY08},
although more sophisticated heuristics exist~\cite{DDFG10,LXY13}.
Then once $\epsilon$ is sufficiently small, 
we fix the sparsity structure by setting small entries in $\vve c$ to zero
and continue optimization until the termination conditions are met.
This final polishing step removes the bias introduced by the regularization.
This process is summarized in \cref{alg:irls}.

\begin{algorithm}[t]
\begin{minipage}{\linewidth}
\begin{algorithm2e}[H]
	\Input{measurements $\vve y$, regularization order $p \ge 0$, 
		basis functions $\lbrace \phi_k\rbrace_{k=1}^n$,
		regularization weight $\alpha \ge 0$, truncation parameter $\tau\ge 0$}
	\Output{
		parameter estimates $\vve c$,
		denoised state $\vve z$
	}
	Initialize $\epsilon^{(0)}\leftarrow 1$,
	$\vve z^{(0)}$ by smoothing~\cref{eq:smooth},
	$\vve c^{(0)}$ by LSOI~\cref{eq:lsopinf} applied to $\vve z^{(0)}$ \;
	\For{$\ell=1,2,3,\ldots$}{
		$\vma W^{(\ell-1)} \leftarrow \diag( |\vve c^{(\ell-1)}|^2 + \epsilon^{(\ell-1)})^{\frac{p}{2}-1}$\;
		Obtain $\vve c^{(\ell)}$ and $\vve z^{(\ell)}$ by one step of SQP applied to \hfill
		$ \displaystyle 
			\min_{\vve c, \vve z} 
				(\vve y - \vve z)^\trans \vma M (\vve y- \vve z) + \alpha \vve c^\trans \vma W^{(\ell-1)} \vve c
			\quad \text{s.t.} \quad  \vma D \vve z = \vma \Phi(\vve z) \vve c
		$\;
		\lIf{optimization converged}{
			$\epsilon^{(\ell)} \leftarrow  \epsilon^{(\ell-1)}/10$
		}
		\lElse{
			$\epsilon^{(\ell)} \leftarrow \epsilon^{(\ell-1)}$
		}
		\lIf{optimization converged and $\epsilon^{(\ell)}< 10^{-8}$}{{\bf break}}
	}
	Fix sparsity structure in $\vve c$ setting $c_i=0$ if $|c_i| \le \tau$\
	and solve
	$ \displaystyle 
		\min_{\vve c, \vve z} 
			(\vve y - \vve z)^\trans \vma M (\vve y- \vve z)
		\quad \text{s.t.} \quad  \vma D \vve z = \vma \Phi(\vve z) \vve c
	$\;
	\vspace*{-1em} 
\end{algorithm2e}
\end{minipage}
\caption{SIDDS with IRLS Regularization Approximation}
\label{alg:irls}
\end{algorithm}

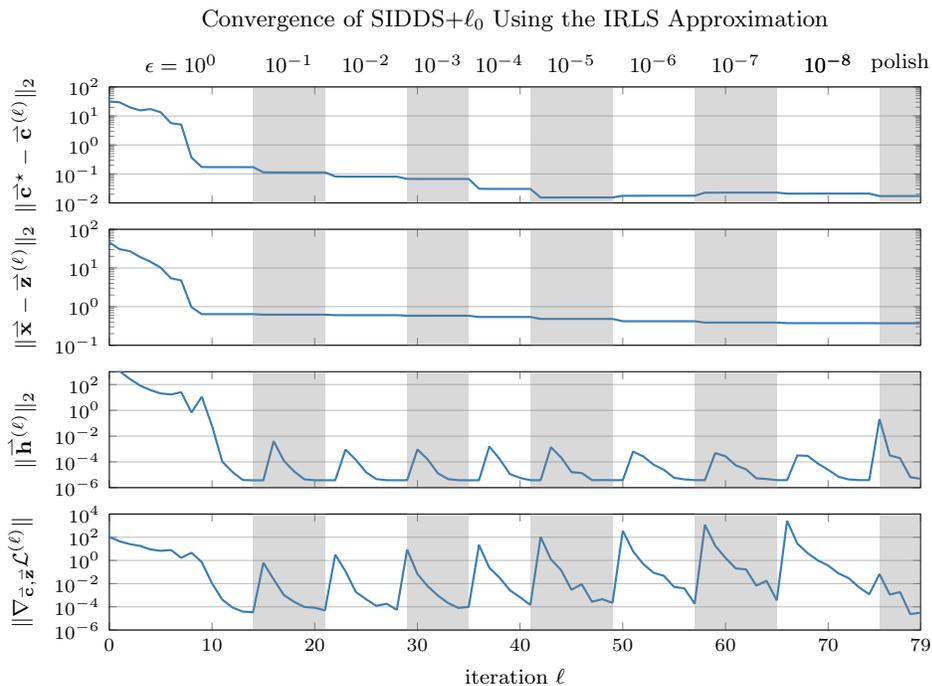
\begin{figure}
\centering

\newcommand{\mydrawpart}{%
	\draw[draw=none, fill = gray, opacity = 0.3] ({axis cs: 14,0} |- {rel axis cs:0,0}) rectangle ({axis cs:21, 0} |- {rel axis cs:(0,1)});
	\draw[draw=none, fill = gray, opacity = 0.3] ({axis cs: 29,0} |- {rel axis cs:0,0}) rectangle ({axis cs:35, 0} |- {rel axis cs:(0,1)});
	\draw[draw=none, fill = gray, opacity = 0.3] ({axis cs: 41,0} |- {rel axis cs:0,0}) rectangle ({axis cs:49, 0} |- {rel axis cs:(0,1)});
	\draw[draw=none, fill = gray, opacity = 0.3] ({axis cs: 57,0} |- {rel axis cs:0,0}) rectangle ({axis cs:65, 0} |- {rel axis cs:(0,1)});
	\draw[draw=none, fill = gray, opacity = 0.3] ({axis cs: 75,0} |- {rel axis cs:0,0}) rectangle ({axis cs:79, 0} |- {rel axis cs:(0,1)});
}

\begin{tikzpicture}
\begin{groupplot}[
	group style = {group size = 1 by 4,
		vertical sep = 1em,
		},
	width = 0.95\linewidth,
	height = 0.24\linewidth,
	ymode = log,
	xmin = 0, xmax = 79,
	xtick = {0,10,20,...,70,79},
	ytickten = {-6,-5, ..., 4},
	ymajorgrids,
	]

	\pgfplotstableread{data/fig_convergence.dat}\conv

	\nextgroupplot[ylabel = $\| \vve c^{\star} - \vve c^{(\ell)}\|_2$, ymin = 1e-2, ymax = 1e2, xticklabels = {,,}, clip = false]
	\addplot[blue, thick] table [x=it, y=C_err] {\conv};
	\mydrawpart
	\node at (axis cs: 7, 1e2) [anchor = south, yshift = 2pt] {$ \epsilon=10^{0}$};  
	\node at (axis cs: 17.5, 1e2) [anchor = south, yshift = 2pt] {$ 10^{-1}$};  
	\node at (axis cs: 25, 1e2) [anchor = south, yshift = 2pt] {$ 10^{-2}$};  
	\node at (axis cs: 32, 1e2) [anchor = south, yshift = 2pt] {$ 10^{-3}$};  
	\node at (axis cs: 38, 1e2) [anchor = south, yshift = 2pt] {$ 10^{-4}$};  
	\node at (axis cs: 45, 1e2) [anchor = south, yshift = 2pt] {$ 10^{-5}$};  
	\node at (axis cs: 53, 1e2) [anchor = south, yshift = 2pt] {$ 10^{-6}$};  
	\node at (axis cs: 61, 1e2) [anchor = south, yshift = 2pt] {$ 10^{-7}$};  
	\node at (axis cs: 70, 1e2) [anchor = south, yshift = 2pt] {$ 10^{-8}$};  
	\node at (axis cs: 70, 1e2) [anchor = south, yshift = 2pt] {$ 10^{-8}$};  
	\node at (axis cs: 77, 1e2) [anchor = south, yshift = 2pt] {polish};  

	
	\nextgroupplot[ylabel = $\| \vve x - \vve z^{(\ell)}\|_2$, ymin = 1e-1, ymax = 1e2, xticklabels = {,,}]
	\addplot[blue, thick] table [x=it, y=X_err] {\conv};
	\mydrawpart
	
	\nextgroupplot[ylabel = $\| \vve h^{(\ell)} \|_2$, ymin = 1e-6, ymax = 1e3, ytickten = {-6, -4, ..., 4}, xticklabels = {,,}]
	\addplot[blue, thick] table [x=it, y=con_norm] {\conv};
	\mydrawpart
	
	\nextgroupplot[ylabel = {$\|\nabla_{\vve c, \vve z} \mathcal{L}^{(\ell)} \| $}, 
		xlabel = iteration $\ell$,
		ymin = 1e-6, ymax = 1e4, ytickten = {-6,-4,..., 4}, ]
	\addplot[blue, thick] table [x=it, y=l_grad_norm] {\conv};
	\mydrawpart
\end{groupplot}
\node (title) at ($(group c1r1.north)+(0,2.5em)$) {\small Convergence of SIDDS+$\ell_0$ Using 
	the IRLS Approximation};
\end{tikzpicture}

\caption[Convergence Plot]{%
Using the IRLS approximation of $R_0$ regularization, 
we rapidly identify the dynamical system to high accuracy.
Stripes denote the different values of $\epsilon$ used as regularization of the IRLS approximation
of the $\ell_0$-norm.
The bottom plot measures the first order optimality of the solution,
where $\mathcal{L}^{(\ell)}$ is the Lagrangian~\cref{eq:lagrangian} at the $\ell$-th iterate.
The data for this problem corresponds to the Lorenz 63 attractor example introduced in \cref{sec:examples:lorenz63}
with additive i.i.d.\ normally distributed noise with standard deviation $\sigma=0.1$ and $m=2000$ measurements.
}
\label{fig:convergence}
\end{figure}

\Cref{fig:convergence} shows how SIDDS proceeds using this IRLS approximation approach.
Starting with large regularization $\epsilon=1$,
the optimizer rapidly identifies a good approximation of the true state,
which then slowly improves over the course of optimization.
Then, as the IRLS regularization term decreases
we slowly obtain better coefficient estimates. 
After each decrease of regularization,
we first see an increase in the Lagrangian gradient norm,
followed by an increase in the constraint mismatch at the subsequent iteration.
However, after only a few more iterations, 
the optimization converges for this particular $\epsilon$
and the value of $\epsilon$ is decreased.

\subsection{Performing SQP Steps\label{sec:algorithm:sqp}}
Using the IRLS convex approximation of the $\ell_p$-norm constraint,
we now seek to solve the corresponding optimization problem,
\begin{equation}\label{eq:opt}
	\begin{split}
	\min_{\vve c, \vve z} \ &  f(\vve c, \vve z) \coloneqq 
		\frac12 (\vve y - \vve z)^\trans \vma M (\vve y - \vve z) 
			+  \frac{\alpha}{2}\vve c^\trans \vma W^{(\ell)} \vve c,\\
	\text{s.t.} \ & \vve h(\vve c, \vve z) \coloneqq \vma D \vve z - \vma \Phi(\vve z) \vve c = \ve 0.
	\end{split}
\end{equation}
Here we use the SQP approach of Liu and Yuan~\cite{LY11}
to ensure convergence without the use of a merit function or filter.
Each step of their algorithm requires the solution of two large-scale subproblems:
a relaxation step and quadratic subproblem.
We show how both can be efficiently solved 
by exploiting the structure of the SIDDS optimization problem.

\subsubsection{Relaxation Step\label{sec:optimization:relax}}
The relaxation step seeks to find a direction $\vve p^{(\ell)}$
that approximately minimizes the error in the linearization of the constraints:
\begin{align}\label{eq:relax2}
	\vve p^{(\ell)} &\approx \argmin_{\vve p} \| \vve h^{(\ell)} + \vma A^{(\ell)} \vve p\|_2^2, \\
	\text{where} \quad
	\vma A^{(\ell)} &\coloneqq \vma A(\vve c^{(\ell)}, \vve z^{(\ell)}),
	\qquad
	\vma A(\vve c, \vve z) \coloneqq 
		\begin{bmatrix} 
				\nabla_{\vve c} \vve h(\vve c, \vve z)&
				\nabla_{\vve z} \vve h(\vve c, \vve z)
		\end{bmatrix},
\end{align}
and $\vve h^{(\ell)} \coloneqq \vve h(\vve c^{(\ell)}, \vve z^{(\ell)})$.
For the convergence analysis of~\cite{LY11} to hold, 
this direction $\vve p^{(\ell)}$ must satisfy two constraints
for small constants $\kappa_1$ and $\kappa_2$
\begin{align}
	\| \vve p^{(\ell)} \| &\le \kappa_1 \| \vma A^{(\ell)\trans} \vve h^{(\ell)}\|, & \kappa_1 &\ge 0; \\
	\| \vve h^{(\ell)} \|^2 - \|\vve h^{(\ell)}\| \|\vve h^{(\ell)} + \vma A^{(\ell)} \vve p^{(\ell)}\| &\ge
		\kappa_2 \|\vma A^{(\ell)\trans} \vve h^{(\ell)}\|^2, &
		\kappa_2 &\in (0,1). 
\end{align}
Due to the scale of $\vma A\in \R^{(dn)\times(dm+dn)}$,
direct solution methods for~\cref{eq:relax2} are impractical.
Instead, we compute the relaxation step by exploiting the structure of $\vma A^{(\ell)}$.
The matrix $\vma A^{(\ell)}$ contains two blocks,
$\vma A^{(\ell)} = \begin{bmatrix} \vma K^{(\ell)} & \vma L^{(\ell)} \end{bmatrix}$,
where 
\begin{align}
	\label{eq:alg:K}
	\vma K^{(\ell)}
		& \coloneqq \vma K(\vve c^{(\ell)}, \vve z^{(\ell)}), &
	\vma K(\vve c, \vve z) 
		& \coloneqq \nabla_{\vve c} \vve h(\vve c, \vve z) = 
		-\vma \Phi(\vve z) ; \\
	\label{eq:alg:L}
	\vma L^{(\ell)}
		& \coloneqq \vma L(\vve c^{(\ell)}, \vve z^{(\ell)}), &
	\vma L(\vve c, \vve z) 
		& \coloneqq \nabla_{\vve z} \vve h(\vve c, \vve z)
		= \vma D - \nabla \vma \Phi(\vve z) \tenmult2 \vve c.
\end{align}
Our approach is to split $\vve p$ into components corresponding to $\vve c$ and $\vve z$:
\begin{equation}\label{eq:split}
		\| \vve h^{(\ell)} + \vma A^{(\ell)} \vve p\|_2^2 
		=
		\| \vve h^{(\ell)} + 
			 \vma K^{(\ell)} \vve p_{\vve c}
			+  \vma L^{(\ell)} \vve p_{\vve z}
		\|_2^2
	\quad \text{where} \quad \vve p = \begin{bmatrix} \vve p_{\vve c} \\ \vve p_{\vve z} \end{bmatrix}
\end{equation}
and then solve for $\vve p_{\vve c}$ and $\vve p_{\vve z}$ separately.

We compute $\vve p_{\vve c}^{(\ell)}$ by setting $\vve p_{\vve z}=\ve 0$ and solving
the overdetermined least squares problem
\begin{align}
	\vve p_{\vve c}^{(\ell)} 
		&\coloneqq \argmin_{\vve p_{\vve c}}\| \vve h^{(\ell)} + \vma K^{(\ell)} \vve p_{\vve c}\|_2^2 \\
		&= \argmin_{\vve p_{\vve c}} \| 
			 \vma \Phi(\vve z^{(\ell)}) \vve p_{\vve c}
				-
			[\vma D \vve z^{(\ell)} - \vma \Phi(\vve z^{(\ell)}) \vve c^{(\ell)}] 
		\|_2^2.
\end{align}
This second statement has the same structure
as LSOI~\cref{eq:lsopinf}
allowing us to restate this in a dense matrix format, removing the Kronecker products:
\begin{equation}\label{eq:relax_sindy}
	\ma P_{\ma C}^{(\ell)} \coloneqq \argmin_{\ma P_{\ma C}} \| 
		\ma \Phi(\ma Z^{(\ell)}) \ma P_{\ma C} - 
		[\ma D \ma Z^{(\ell)} - \ma \Phi(\ma Z^{(\ell)}) \ma C^{(\ell)}]
		\|_\fro^2.
\end{equation}
This allows more efficient solution via a QR factorization
and we set $\vve p_{\vve c}^{(\ell)} = \vectorize(\ma P_{\ma C}^{(\ell)})$.

Next, we solve for $\vve p_{\vve z}$ holding $\vve p_{\vve c}$ constant:
\begin{equation}
	\vve p_{\vve z}^{(\ell)} \approx \argmin_{\vve p_{\vve z}}
		\| \vve h^{(\ell)} + \vma K^{(\ell)} \vve p_{\vve c}^{(\ell)}
			+ \vma L^{(\ell)} \vve p_{\vve z}\|_2^2.
\end{equation}
Although $\vma L^{(\ell)}$ is square, 
it is structurally rank deficient.
The matrix $\vma L^{(\ell)}$ encodes the discretized constraint for 
$\op E^{\delta (j-1)}(\ve z_1, \vve c) = \ve z_j$ for $j=2, \ldots, m$;
this continuous constraint only provides $d(m-1)$ constraints whereas $\vma L^{(\ell)}$ encodes $dm$ constraints.
Thus when approximating $\vve p_{\vve z}$ we include a small amount of Tikhonov regularization and solve via the normal equations:
\begin{equation}
	\vve p_{\vve z}^{(\ell)} \coloneqq
		 -[\vma L^{(\ell)\trans} \vma L^{(\ell)} + \beta \ma I]^{-1}
		\vma L^{(\ell)\trans} [\vve h^{(\ell)} + \vma K^{(\ell)} \vve p_{\vve c}^{(\ell)}];
\end{equation}
in our numerical experiments we take $\beta = 10^{-6}$.
Since $\vma L^{(\ell)}\in \R^{(dm)\times (dm)}$ has small bandwidth $dq \ll dm$,
we efficiently apply its inverse using a sparse LU factorization.

\subsubsection{Solution of Quadratic Subproblem\label{sec:optimization:qp}}
The other expensive component of the Liu and Yuan SQP algorithm 
is solving the relaxed quadratic program
\begin{equation}\label{eq:subqp}
	\begin{split}
		\vve d^{(\ell)} = 
		\argmin_{\vve d}  & \ {\vve g^{(\ell)}}^\trans \vve d + \frac12 \vve d^\trans \vma B^{(\ell)} \vve d \\
		\quad \text{s.t.} & \  \vma A^{(\ell)} \vve d = \vma A^{(\ell)} \vve p^{(\ell)}
	\end{split}
\end{equation}
where $\vma B^{(\ell)}$ is an approximation of the Lagrangian Hessian
and $\vve g^{(\ell)}$ is the gradient
\begin{equation}
	\vve g^{(\ell)} \coloneqq \vve g(\vve c^{(\ell)}, \vve z^{(\ell)})
	\qquad
	\vve g(\vve c, \vve z) 
	\coloneqq
		\begin{bmatrix} 
			\nabla_{\vve c} f(\vve c, \vve z) &
			\nabla_{\vve z} f(\vve c, \vve z)
		\end{bmatrix}.
\end{equation}
Here we introduce a sparse approximation of the $\vma B^{(\ell)}$
and show how to efficiently solve~\cref{eq:subqp}
by direct solution of the stabilized KKT system with a preconditioned MINRES iteration.

Ideally we would use the exact Hessian of the Lagrangian for $\vve B^{(\ell)}$
\begin{equation}\label{eq:lagrangian}
	\mathcal{L}(\vve c, \vve z, \vve w) \coloneqq f(\vve c, \vve z) + \vve w^\trans \vve h(\vve c, \vve z),
\end{equation}
where $\vve w$ are the Lagrange multipliers.
However, this choice is impractical as $\nabla_{\vve c, \vve z}^2\mathcal{L}$ will generally be dense.
Note the contribution from the Lagrange multiplier term is
\begin{equation}
	\vve w^\trans \vve h(\vve c, \vve z) = \vve w^\trans \vma D_p \vve z - 
	\vve w^\trans \vma \Phi(\vve z) \vve c;
\end{equation}
the Hessian of this second term $\vve w^\trans \vma \Phi(\vve z) \vve c$
will be dense except for special combinations of $\vve c$ and basis vectors $\phi_k$.
Instead we neglect this term entirely 
and approximate the Hessian of the Lagrangian by the Hessian of the objective:
\begin{equation}
	\vma B^{(\ell)} \coloneqq \nabla_{\vve c, \vve z}^2 f(\vve c^{(\ell)}, \vve z^{(\ell)}) 
		= 
	\begin{bmatrix}
		\alpha \vma W^{(\ell)} & \ma 0 \\ \ma 0 & \vma M
	\end{bmatrix}.
\end{equation}
By construction, this is a positive semidefinite matrix.

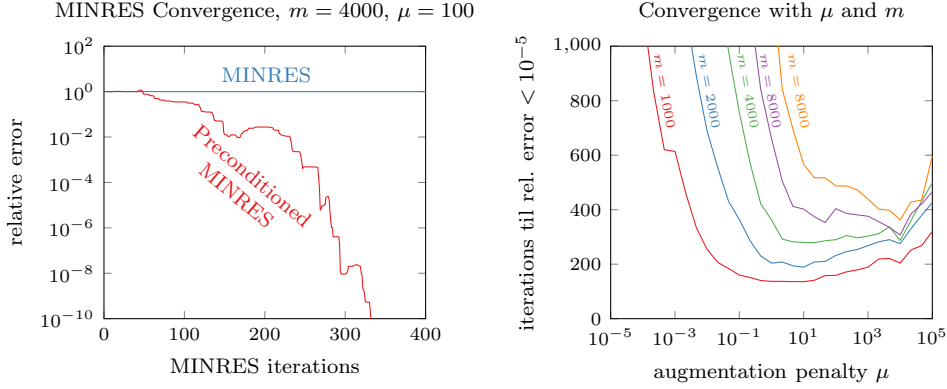
\begin{figure}
\centering

\begin{tikzpicture}
\begin{groupplot}[
	width = 0.45\linewidth,
	height = 0.4\linewidth,
	group style = {
		group size = 2 by 1,
		horizontal sep = 7em,
	},
]
\nextgroupplot[
	ymin = 1e-10,
	ymax = 1e2,
	ytickten = {-14, -12, ..., 2},
	ymode = log,
	xmin = 0, 
	xmax = 400,
	xtick = {0, 100, ..., 500},
	xlabel = MINRES iterations,
	ylabel = relative error, 
	title = {MINRES Convergence, $m=4000$, $\mu=100$},
]
	\addplot[blue] table [x expr = \lineno + 1, y = dist] {data/fig_minres_converge_M_4000_mu_1.000000e+02.dat}
		node [pos = 0.1, anchor = south] {\footnotesize MINRES}; 
	
	\addplot[red] table [x expr = \lineno + 1, y = dist] {data/fig_minres_converge_pre_M_4000_mu_1.000000e+02.dat}
		node [pos = 0.6, anchor = north, rotate=-40, yshift=-12pt] {\begin{tabular}{c}Preconditioned \\ MINRES\end{tabular}}; 
	

\nextgroupplot[
	xlabel = augmentation penalty $\mu$,
	ylabel = iterations til rel. error $<10^{-5}$,
	xmode = log,
	xmin = 1e-5,
	xmax = 1e5,
	xtickten = {-5, -3, -1, 1, 3, 5},
	ytick = {0, 200, ..., 2000},
	ymin = 0,
	ymax = 1000,
	title = {Convergence with $\mu$ and $m$},
]
	\addplot table [x = mu, y = it_1.000000e-05] {data/fig_minres_mu_M_1000.dat};
	\draw[red] (1e-4,1000) node [anchor = south west, rotate = -80] {\tiny $m=1000$};
	
	\addplot table [x = mu, y = it_1.000000e-05] {data/fig_minres_mu_M_2000.dat}; 
	\draw[blue] (2.5e-3,1000) node [anchor = south west, rotate = -80] {\tiny $m=2000$};
	
	\addplot table [x = mu, y = it_1.000000e-05] {data/fig_minres_mu_M_4000.dat};
	\draw[green] (4e-2,1000) node [anchor = south west, rotate = -80] {\tiny $m=4000$};
	
	\addplot table [x = mu, y = it_1.000000e-05] {data/fig_minres_mu_M_8000.dat};
	\draw[purple] (0.2,1000) node [anchor = south west, rotate = -80] {\tiny $m=8000$};

	\addplot table [x = mu, y = it_1.000000e-05] {data/fig_minres_mu_M_16000.dat};
	\draw[orange] (1.5,1000) node [anchor = south west, rotate = -80] {\tiny $m=8000$};

\end{groupplot}
\end{tikzpicture}

\caption[MINRES]{%
	The preconditioner $\op M$ accelerates the convergence of MINRES on the KKT system~\cref{eq:KKT}.
	On the left we show the convergence the MINRES iteration
	by showing the 2-norm difference between each iterate's estimate of the solution
	and the converged value.
	On the right, we see the number of MINRES iterations required varies as a function of 
	augmentation penalty $\mu$ and problem dimension $m$.
	Test data comes from the Van der Pol oscillator 
	described in \cref{sec:examples:vanderpol}.
}
\label{fig:minres}
\end{figure}

Next, we seek to solve the KKT system
\begin{align}
	\begin{bmatrix}
		\vma B^{(\ell)} & \vma A^{(\ell)\trans} \\
		\vma A^{(\ell)} & \ma 0
	\end{bmatrix}
	\begin{bmatrix}
		\vve d \\
		\vve w
	\end{bmatrix}
	=
	\begin{bmatrix}
		-\vve g^{(\ell)} \\
		\vma A^{(\ell)}\vve p^{(\ell)}
	\end{bmatrix}.
\end{align}
Since $\vma L^{(\ell)}$ is rank deficient,
$\vma A^{(\ell)}$ does not have full row rank
causing the KKT system to be singular.
To correct this, we use a stabilization procedure following Wright~\cite[eq.~(4.2)]{Wri05}
where we add multiple of the identity to the $(1,1)$ and $(2,2)$ blocks,
\begin{align}\label{eq:KKT}
	\begin{bmatrix}
		\vma B^{(\ell)} + \zeta \ma I & \vma A^{(\ell)\trans} \\
		\vma A^{(\ell)} & -\gamma^{(\ell)} \ma I
	\end{bmatrix}
	\begin{bmatrix}
		\vve d^{(\ell)} \\
		\vve w^{(\ell)}
	\end{bmatrix}
	=
	\begin{bmatrix}
		-\vve g^{(\ell)} \\
		\vma A^{(\ell)}\vve p^{(\ell)}
	\end{bmatrix}.
\end{align}
In our implementation we fix $\zeta = 10^{-4}$
and scale $\gamma$ based on the previous iteration
\begin{equation}
	\gamma^{(\ell)} = 10^{-4} \left( \|\vve g^{(\ell-1)} + \vve w^{(\ell-1)} \vma A^{(\ell-1)} \|_1 + \| \vve h^{(\ell-1)}\|_1 \right).
\end{equation}
To solve~\cref{eq:KKT},
we use a preconditioned MINRES iteration following Alger, et al.~\cite{AVBG17}.
This preconditioner is derived from a block-diagonal approximation
of the KKT system corresponding the augmented objective,
\begin{equation}
	f_{\mu}(\vve c, \vve z) \coloneqq f(\vve c, \vve z) + \mu \| \vve h(\vve c, \vve w)\|_2^2
\end{equation}
for some penalty $\mu > 0$.
In our setting, this preconditioner is
\begin{equation}
\op{M}^{(\ell)} =
\begin{bmatrix}
	\alpha \vma W + \zeta \ma I + \mu \vma K^{(\ell)\trans} \vma K^{(\ell)} & & \\
	& \vma M + \zeta \ma I + \mu \vma L^{(\ell)\trans} \vma L^{(\ell)} &  \\
	& & \frac{1}{\mu} \ma I 
\end{bmatrix}.
\end{equation}
We can efficiently apply ${\op{M}^{(\ell)}}^{-1}$ blockwise.
The first block is small, $dn$, so we explicitly compute its inverse.
The second block is large, $dm$, but low bandwidth $dq \ll dm$,
so we apply its inverse by precomputing a sparse factorization;
in our experiments an LU factorization computed using SuperLU~\cite{superlu99}.
As illustrated in \cref{fig:minres}, 
this preconditioner enables rapid solution of the KKT system.
In our setting,
the performance of this preconditioner depends on the augmentation penalty $\mu$.
These experiments show that a value $\mu \in [10,1000]$ enables fast convergence 
for a variety of problem dimensions;
our experiments take $\mu =100$.

\subsection{Multiple Trajectories\label{sec:optimization:traj}}
In some situations it is necessary incorporate data from multiple trajectories
to ensure the operator inference problem is well posed~\cite{WX19}.
For example, a single trajectory might not sufficiently explore the state-space
to enable an accurate estimate of the parameters $\vve c$.
Fortunately, SIDDS can be easily modified to 
accommodate this situation.
Suppose we have trajectories $\lbrace \vve y_{(i)}\rbrace_{i=1}^N$,
we then seek to solve an extension of~\cref{eq:opt}
\begin{equation}\label{eq:opt_mult}
	\begin{split}
	\min_{\vve c, \vve z_{(1)}, \ldots, \vve z_{(N)}} \ & 
		\sum_{i=1}^N \frac12 [\vve y_{(i)} - \vve z_{(i)}]^\trans \vma M_{(i)} [\vve y_{(i)} - \vve z_{(i)}] 
		+
		\frac{\alpha}{2}  \vve c^\trans \vma W \vve c, \\
	\text{s.t.} \ & \vma D_p\vve z_{(i)} - \vma \Phi(\vve z_{(i)}) \vve c = \ve 0, \quad i=1,\ldots, N.
	\end{split}
\end{equation}
Structurally, this optimization problem is similar to~\cref{eq:opt}
enabling us to use the same techniques to efficiently solve this problem.
We build the constraint derivative 
using components defined analogously to~\cref{eq:alg:K} and~\cref{eq:alg:L},
\begin{align}
	\vma K_{(i)}^{(\ell)} &\coloneqq -\vma \Phi(\vve z_{(i)}^{(\ell)}), & 
	\vma L_{(i)}^{(\ell)} &\coloneqq \vma D_p - 
		[ \nabla \vma \Phi(\vve z_{(i)}^{(\ell)}) \tenmult2 \vve c^{(\ell)}],
\end{align}
forming the full constraint derivative,
\begin{equation}
	\vma A^{(\ell)} \coloneqq \begin{bmatrix}
		\vma K_{(1)}^{(\ell)} & \vma L_{(1)}^{(\ell)}  \\
		\vdots && \ddots \\
		\undermat{\vma K_{\square}^{(\ell)}}{\vma K_{(N)}^{(\ell)}} &
		\undermat{\vma L_{\square}^{(\ell)}}{ \phantom{\vma L_{(1)}^{(\ell)}} %
			& \phantom{\vma L_{(2)}^{(\ell)}} & 
		\vma L_{(N)}^{(\ell)}} 
	\end{bmatrix}
	\vphantom{\begin{bmatrix} \vma L \\ \vma L \\ \vma L \\ \vma L  \\ \vma L \end{bmatrix}}.
\end{equation}
These two matrices $\vma K_\square^{(\ell)}$ and $\vma L_\square^{(\ell)}$
with a dense rectangular block and a sparse, low-bandwidth square matrix
can be used analogously in place of $\vma K^{(\ell)}$ and $\vma L^{(\ell)}$ in the preceding analysis
for the relaxation step and the quadratic subproblem.
This enables an efficient solution to SIDDS with multiple trajectory data.


\section{Numerical Experiments\label{sec:examples}}
In this section we provide a few numerical experiments illustrating the performance of SIDDS
and comparing against existing methods 
on three representative test problems.
In each case, we use basis functions $\lbrace \phi_k\rbrace_k$
corresponding to a degree-$p$ total-degree monomial basis,
\begin{equation}
	\mathcal{P}_d^p
	\coloneqq \lbrace \phi_{\ve \alpha}(\ve x) \rbrace_{|\ve \alpha| \le p}, 
	\qquad
	\phi_{\ve \alpha}(\ve x) = \prod_{i=1}^d x_i^{\alpha_i},
	\qquad 
	|\ve \alpha| \coloneqq \sum_{i=1}^d \alpha_i
\end{equation}
where $\ve \alpha \in \mathbb{Z}_+^d$ is a multi-index over $d$ nonnegative integers.
In all our experiments we choose 
the basis $\mathcal{P}_d^p$ with the smallest $p$
containing the example differential equation.

\subsection{Test Problems\label{sec:examples:test}}
Here we consider three common test problems
the Duffing oscillator, Lorenz 63 attractor, and the Van der Pol Oscillator;
see, e.g.,~\cite{BPK16,CPD21,KBK20x}.
We focus on these low-dimensional problems
were we can perform Monte Carlo tests 
to evaluate performance over multiple realizations of noise.
Unless otherwise mentioned,
we use the sample rate $\delta=10^{-2}$.
To generate measurements $\vve y$,
we first evolve the corresponding differential equation using SciPy's 
\verb|solve_ivp| with the \verb|DOP853| integrator to generate $\vve y$.
Then to generate noisy measurements,
we sample $\vve n \sim \mathcal{N}(\vve 0, \vma \Sigma)$
and form $\vve y = \vve x + \vve n$.

\subsubsection{Duffing Oscillator}
The Duffing oscillator models a nonlinear pendulum;
here we considered a damped variant where
\begin{equation}
	\left\lbrace
	\begin{aligned}
		\dot{x}_1 &= x_2, & 
		x_1(0) &= -2, \\
		\dot{x}_2 &= -0.1 x_2 -  x_1 - 5 x_1^3, \quad &
		x_2(0) &= -2.
	\end{aligned}
	\right.
\end{equation}
In our experiments,
we typically use $m=1000$ measurements of this system.

\subsubsection{Lorenz 63 Attractor\label{sec:examples:lorenz63}}
The Lorenz 63 attractor is a chaotic system 
developed initially from models of atmospheric convection:
\begin{equation}
	\left\lbrace
	\begin{aligned}
	\begin{split}
		\dot{x}_1 &= 10(x_2 - x_1),  &
		x_1(0) &= -8, \\
		\dot{x}_2 &= x_1(28 - x_3) - x_2, \quad &
		x_2(0) &= 7, \\
		\dot{x}_3 &= x_1 x_2 -\tfrac83 x_3, &
		x_3(0) &= -28. 
	\end{split}
	\end{aligned}
	\right.
\end{equation}
In our experiments, we typically use $m=2000$ measurements 
of this system.

\subsubsection{Van der Pol Oscillator\label{sec:examples:vanderpol}}
The Van der Pol oscillator is a nonlinear ODE 
with a non-trivial limit cycle
\begin{equation}
	\left\lbrace
	\begin{aligned}
		\dot{x}_1 &= x_2, & 
		x_1(0) &= 0, \\
		\dot{x}_2 &= 2 x_2(1 - x_1^2) - x_1, \quad &
		x_2(0) &= 1.
	\end{aligned}
	\right.
\end{equation}
In our experiments,
we typically use $m=1000$ measurements of this system.

\subsection{Tuning ODE Integration\label{sec:examples:integration}}
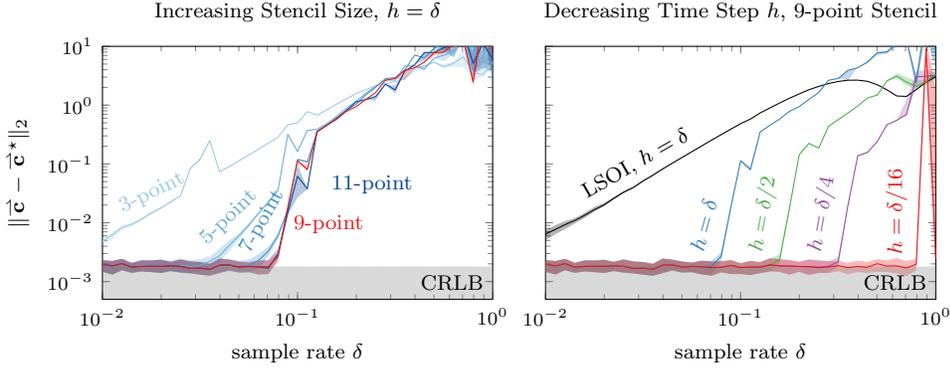
\begin{figure}
\centering

\begin{tikzpicture}
\begin{groupplot}[
	xmode = log,
	ymode = log,
	width = 0.52\linewidth,
	height = 0.38\linewidth,
	ymin = 5e-4,
	ymax = 1e1,
	ytickten = {-3,-2,...,1},
	xlabel = {sample rate $\delta$},
	ylabel = {$\| \vve c - \vve c^\star\|_2$},
	group style = {group size = 2 by 1,
		horizontal sep = 2em,
		},
]

\pgfplotstableread{data/fig_integration.dat}\crlb

\nextgroupplot[title={Increasing Stencil Size, $h=\delta$}]

\pgfplotstableread{data/fig_sample_rate_lide_points_3_oversample_1.dat}\lidethree
\addplot[Blues7C] table [x=dt, y=F50] {\lidethree}
	node [pos=0.05, anchor = south west, rotate=26, yshift=0pt, xshift=-2pt] {3-point};
\addplot[draw=none, name path=lidethree_lb] table [x=dt, y=F25] {\lidethree};
\addplot[draw=none, name path=lidethree_ub] table [x=dt, y=F75] {\lidethree};
\addplot[Blues7C, opacity=0.3] fill between [of=lidethree_lb and lidethree_ub];

\pgfplotstableread{data/fig_sample_rate_lide_points_5_oversample_1.dat}\lidefive
\addplot[Blues7D] table [x=dt, y=F50] {\lidefive}
	node [pos=.1, anchor = south west, rotate=44, yshift=0pt] {5-point};
\addplot[draw=none, name path=lidefive_lb] table [x=dt, y=F25] {\lidefive};
\addplot[draw=none, name path=lidefive_ub] table [x=dt, y=F75] {\lidefive};
\addplot[Blues7D, opacity=0.3] fill between [of=lidefive_lb and lidefive_ub];

\pgfplotstableread{data/fig_sample_rate_lide_points_7_oversample_1.dat}\lideseven
\addplot[Blues7E] table [x=dt, y=F50] {\lideseven}
	node [pos=0.13, anchor = south west, rotate=60, yshift=-3pt, xshift=0pt] {7-point};
\addplot[draw=none, name path=lideseven_lb] table [x=dt, y=F25] {\lideseven};
\addplot[draw=none, name path=lideseven_ub] table [x=dt, y=F75] {\lideseven};
\addplot[Blues7E, opacity=0.3] fill between [of=lideseven_lb and lideseven_ub];

\pgfplotstableread{data/fig_sample_rate_lide_points_11_oversample_1.dat}\lideeleven
\addplot[Blues7G] table [x=dt, y=F50] {\lideeleven}
	node [pos=.3, anchor = south west, rotate=0, xshift=12pt] {11-point};
\addplot[draw=none, name path=lideeleven_lb] table [x=dt, y=F25] {\lideeleven};
\addplot[draw=none, name path=lideeleven_ub] table [x=dt, y=F75] {\lideeleven};
\addplot[Blues7G, opacity=0.3] fill between [of=lideeleven_lb and lideeleven_ub];


\pgfplotstableread{data/fig_sample_rate_lide_points_9_oversample_1.dat}\lidenine
\addplot[red] table [x=dt, y=F50] {\lidenine}
	node [pos=.22, anchor = west, rotate=0, yshift=0pt] {9-point};
\addplot[draw=none, name path=lidenine_lb] table [x=dt, y=F25] {\lidenine};
\addplot[draw=none, name path=lidenine_ub] table [x=dt, y=F75] {\lidenine};
\addplot[red, opacity=0.3] fill between [of=lidenine_lb and lidenine_ub];

\addplot[draw=none,  name path = crb] 
	table [x=dt, y expr = 1e-3*\thisrow{crlb_median}]{\crlb};
\addplot[draw=none, name path = crblb, domain=1e-2:1] {5e-4}
	node [pos=1, anchor = south east, rotate=0] {CRLB};
\addplot[gray, opacity =0.3] fill between [of=crblb and crb];

\nextgroupplot[title={Decreasing Time Step $h$, 9-point Stencil}, ylabel = {}, yticklabels = {,,}]

\pgfplotstableread{data/fig_sample_rate_lsopinf_points_3_oversample_1.dat}\lsopinf
\addplot[black] table [x=dt, y=F50] {\lsopinf}
	node [pos=0.1, anchor = south west, rotate=30, yshift=0pt] {LSOI, $h=\delta$};
\addplot[draw=none, name path=lsopinf_lb] table [x=dt, y=F25] {\lsopinf};
\addplot[draw=none, name path=lsopinf_ub] table [x=dt, y=F75] {\lsopinf};
\addplot[black, opacity=0.3] fill between [of=lsopinf_lb and lsopinf_ub];

\pgfplotstableread{data/fig_sample_rate_lide_points_9_oversample_1.dat}\lideA
\addplot[blue] table [x=dt, y=F50] {\lideA};
\node at (axis cs:7e-2, 2e-3) [anchor=south west, rotate = 75, blue] {$h=\delta$};
\addplot[draw=none, name path=lideA_lb] table [x=dt, y=F25] {\lideA};
\addplot[draw=none, name path=lideA_ub] table [x=dt, y=F75] {\lideA};
\addplot[blue, opacity=0.3] fill between [of=lideA_lb and lideA_ub];

\pgfplotstableread{data/fig_sample_rate_lide_points_9_oversample_2.dat}\lideB
\addplot[green] table [x=dt, y=F50] {\lideB};
\node at (axis cs:1.5e-1, 2e-3) [anchor=south west, rotate = 80, green] {$h=\delta/2$};
\addplot[draw=none, name path=lideB_lb] table [x=dt, y=F25] {\lideB};
\addplot[draw=none, name path=lideB_ub] table [x=dt, y=F75] {\lideB};
\addplot[green, opacity=0.3] fill between [of=lideB_lb and lideB_ub];

\pgfplotstableread{data/fig_sample_rate_lide_points_9_oversample_4.dat}\lideD
\addplot[purple] table [x=dt, y=F50] {\lideD};
\node at (axis cs:3e-1, 2e-3) [anchor=south west, rotate = 80, purple] {$h=\delta/4$};
\addplot[draw=none, name path=lideD_lb] table [x=dt, y=F25] {\lideD};
\addplot[draw=none, name path=lideD_ub] table [x=dt, y=F75] {\lideD};
\addplot[purple, opacity=0.3] fill between [of=lideD_lb and lideD_ub];

\pgfplotstableread{data/fig_sample_rate_lide_points_9_oversample_16.dat}\lide
\addplot[red] table [x=dt, y=F50] {\lide};
\node at (axis cs:8e-1, 2e-3) [anchor=south west, rotate = 90, red] {$h=\delta/16$};
\addplot[draw=none, name path=lide_lb] table [x=dt, y=F25] {\lide};
\addplot[draw=none, name path=lide_ub] table [x=dt, y=F75] {\lide};
\addplot[red, opacity=0.3] fill between [of=lide_lb and lide_ub];

\addplot[draw=none,  name path = crb] 
	table [x=dt, y expr = 1e-3*\thisrow{crlb_median}]{\crlb};
\addplot[draw=none, name path = crblb, domain=1e-2:1] {5e-4}
	node [pos=1, anchor = south east, rotate=0] {CRLB};
\addplot[gray, opacity =0.3] fill between [of=crblb and crb];

\end{groupplot}
\end{tikzpicture}

\caption[Integration]{
	Decreasing the discretization error in the ODE constraint of SIDDS
	improves the ability to recover a dynamical system.
	We use two techniques:
	increasing the order of the derivative approximation (left)
	and decreasing the time-step (right).
	In this example we collect $m=1000$ measurements of the Van der Pol system 
	starting from a point on its limit cycle,
	adding i.i.d.\ normal noise with standard deviation $\sigma = 10^{-3}$.
	We estimate $\vve c$ using SIDDS with $\vma M = \vma I$ on the left
	and an expanded $\vma M^\uparrow$ on the right.
	The shaded region encloses the range between the 25th to 75th percentiles
	from 100 trials; the solid line indicates the median.
	The gray shaded region indicates
	the lower bound on performance 
	given by the Cram\'er-Rao Lower Bound (CRLB).
}
\label{fig:integration}
\end{figure}

Both SIDDS and LSOI introduce a bias in the coefficients $\vve c$
proportional to the discretization error of the ODE.
With SIDDS though we have two techniques we can employ 
to reduce this error:
increasing the accuracy of the derivative approximation by using a larger 
finite difference stencil
and decreasing the time step in the ODE constraint.
Both techniques are illustrated in \cref{fig:integration}
with the Van der Pol example
where we decrease the sample rate to increase the discretization error.
In this example we also initialize the system at a point on its limit cycle
so that the CRLB remains approximately constant irrespective of sample rate.

Decreasing discretization error by using a higher order
finite difference stencil is successful, 
but yields diminishing returns beyond a 9-point stencil
as this example illustrates.
Moreover, larger stencils increase computational cost because
increasing the stencil width increases the bandwidth of $\vma L$~\cref{eq:alg:L} 
and consequently the cost of its LU factorization.
We use this 9-point stencil in all remaining experiments.

Another way to decrease the discretization error
is to decrease the time step.
With SIDDS we can set the ODE integration time step $h$
to be any positive integer fraction of the sample rate $\delta$.
For example, if we want to integrate with time step $h = \delta/2$,
we introduce a zero vector between every measurement $\ve y_j$
and a corresponding zero block in the weight matrix $\vma M$;
assuming $\vma M=\ma I$,
\begin{align}\label{eq:oversample}
	\vve y^{\uparrow} &= 
	\begin{bmatrix} 
	\ve y_1^\trans & \ve 0 & \ve y_2^\trans & \ve 0 & \cdots & \ve y_m^\trans 
	\end{bmatrix}^\trans, &
	\vma M^{\uparrow} &= \diag(\ma I_d, \ma 0_d, \ma I_d, \ma 0_d, \cdots, \ma I_d),
\end{align}
where $\ma I_d$ and $\ma 0_d$ are the $d\times d$ identity and zeros matrices.
This new $\vve y^\uparrow$ has (effectively) half the sample rate, $\delta^\uparrow = h = \delta/2$,
and the zeros in $\vma M^\uparrow$ ensure the values added between measurements
do not affect the objective.
As \cref{fig:integration} shows,
solving SIDDS with these expanded quantities $\vve y^\uparrow$ and $\vma M^\uparrow$
allows the accurate identification
even with slow sample rates $\delta$---%
something not possible with LSOI.
Although useful when sample rate is slow relative to the dynamics,
decreasing the time step substantially increases the cost solving the SIDDS problem.
As our remaining examples are not in this regime, 
we use the same time step as the sample rate ($h=\delta$) in the rest of our experiments.


\subsection{Correlated Noise\label{sec:examples:correlated}}
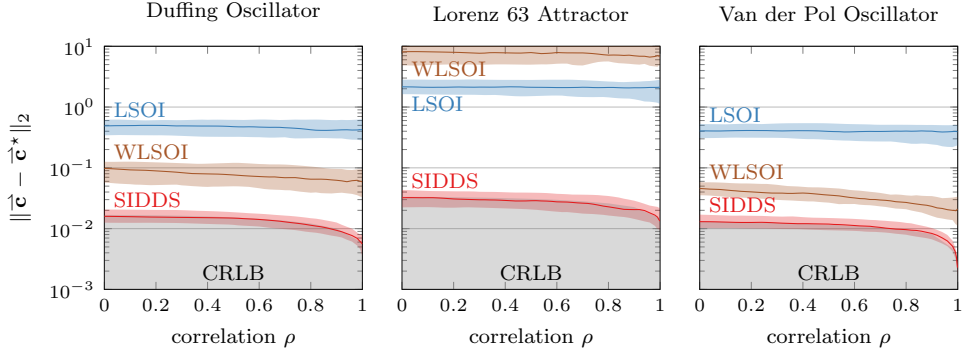
\begin{figure}

\centering
\begin{tikzpicture}
\begin{groupplot}[
	group style = {group size = 3 by 1,
		horizontal sep = 1.5em,
		vertical sep = 2em,
		},
	xmode = linear,
	ymode = log,
	width = 0.385\linewidth,
	height = 0.37\linewidth,
	ytickten = {-7,-6, ..., 5},
	ymajorgrids,
	xlabel = correlation $\rho$,
	xmin = 0,
	xmax = 1,
	ymin = 1e-3,
	ymax = 1e1
]
\nextgroupplot[title = Duffing Oscillator, ylabel = $\| \vve c - \vve c^\star\|_2$ ]

\pgfplotstableread{data/fig_correlated_duffing_lsopinf.dat}\lsopinf
\addplot[blue] table [x=correlation, y=F50] {\lsopinf}
	node [pos=0, anchor = south west, rotate=0] {LSOI};
\addplot[draw=none, name path=lsopinf_lb] table [x=correlation, y=F25] {\lsopinf};
\addplot[draw=none, name path=lsopinf_ub] table [x=correlation, y=F75] {\lsopinf};
\addplot[blue, opacity=0.3] fill between [of=lsopinf_lb and lsopinf_ub];

\pgfplotstableread{data/fig_correlated_duffing_wlsoi.dat}\wlsoi
\addplot[brown] table [x=correlation, y=F50] {\wlsoi}
	node [pos=0, anchor = south west , rotate=0] {WLSOI};
\addplot[draw=none, name path=wlsoi_lb] table [x=correlation, y=F25] {\wlsoi};
\addplot[draw=none, name path=wlsoi_ub] table [x=correlation, y=F75] {\wlsoi};
\addplot[brown, opacity=0.3] fill between [of=wlsoi_lb and wlsoi_ub];

\pgfplotstableread{data/fig_correlated_duffing_lide.dat}\lide
\addplot[red] table [x=correlation, y=F50] {\lide}
	node [pos=0, anchor = south west, rotate=0] {SIDDS};
\addplot[draw=none, name path=lide_lb] table [x=correlation, y=F25] {\lide};
\addplot[draw=none, name path=lide_ub] table [x=correlation, y=F75] {\lide};
\addplot[red, opacity=0.3] fill between [of=lide_lb and lide_ub];

\pgfplotstableread{data/fig_correlated_duffing_crlb.dat}\crlb
\addplot[draw=none, domain = 1e-3:1, name path = crb] table[x=correlation, y expr=\thisrow{crlb_median}*1e-2] 
	{\crlb};
\addplot[draw=none, name path = crblb, domain=0:1] {1e-3}
	node [pos=0.5, anchor = south, rotate=0] {CRLB};
\addplot[gray, opacity =0.3] fill between [of=crblb and crb];

\nextgroupplot[title = Lorenz 63 Attractor, yticklabels = {,,} ]
\pgfplotstableread{data/fig_correlated_lorenz63_lsopinf.dat}\lsopinf
\addplot[blue] table [x=correlation, y=F50] {\lsopinf}
	node [pos=0, anchor = north west, rotate=0] {LSOI};
\addplot[draw=none, name path=lsopinf_lb] table [x=correlation, y=F25] {\lsopinf};
\addplot[draw=none, name path=lsopinf_ub] table [x=correlation, y=F75] {\lsopinf};
\addplot[blue, opacity=0.3] fill between [of=lsopinf_lb and lsopinf_ub];

\pgfplotstableread{data/fig_correlated_lorenz63_wlsoi.dat}\wlsoi
\addplot[brown] table [x=correlation, y=F50] {\wlsoi}
	node [pos=0, anchor = north west, rotate=0] {WLSOI};
\addplot[draw=none, name path=wlsoi_lb] table [x=correlation, y=F25] {\wlsoi};
\addplot[draw=none, name path=wlsoi_ub] table [x=correlation, y=F75] {\wlsoi};
\addplot[brown, opacity=0.3] fill between [of=wlsoi_lb and wlsoi_ub];

\pgfplotstableread{data/fig_correlated_lorenz63_lide.dat}\lide
\addplot[red] table [x=correlation, y=F50] {\lide}
	node [pos=0, anchor = south west, rotate=0] {SIDDS};
\addplot[draw=none, name path=lide_lb] table [x=correlation, y=F25] {\lide};
\addplot[draw=none, name path=lide_ub] table [x=correlation, y=F75] {\lide};
\addplot[red, opacity=0.3] fill between [of=lide_lb and lide_ub];

\pgfplotstableread{data/fig_correlated_lorenz63_crlb.dat}\crlb
\addplot[draw=none, domain = 1e-3:1, name path = crb] table[x=correlation, y expr=\thisrow{crlb_median}*1e-2] 
	{\crlb};
\addplot[draw=none, name path = crblb, domain=0:1] {1e-3}
	node [pos=0.5, anchor = south, rotate=0] {CRLB};
\addplot[gray, opacity =0.3] fill between [of=crblb and crb];

\nextgroupplot[title = Van der Pol Oscillator, yticklabels = {,,}]
\pgfplotstableread{data/fig_correlated_vanderpol_lsopinf.dat}\lsopinf
\addplot[blue] table [x=correlation, y=F50] {\lsopinf}
	node [pos=0, anchor = south west, rotate=0] {LSOI};
\addplot[draw=none, name path=lsopinf_lb] table [x=correlation, y=F25] {\lsopinf};
\addplot[draw=none, name path=lsopinf_ub] table [x=correlation, y=F75] {\lsopinf};
\addplot[blue, opacity=0.3] fill between [of=lsopinf_lb and lsopinf_ub];

\pgfplotstableread{data/fig_correlated_vanderpol_wlsoi.dat}\wlsoi
\addplot[brown] table [x=correlation, y=F50] {\wlsoi}
	node [pos=0, anchor = south west, rotate=0] {WLSOI};
\addplot[draw=none, name path=wlsoi_lb] table [x=correlation, y=F25] {\wlsoi};
\addplot[draw=none, name path=wlsoi_ub] table [x=correlation, y=F75] {\wlsoi};
\addplot[brown, opacity=0.3] fill between [of=wlsoi_lb and wlsoi_ub];

\pgfplotstableread{data/fig_correlated_vanderpol_lide.dat}\lide
\addplot[red] table [x=correlation, y=F50] {\lide}
	node [pos=0, anchor = south west, rotate=0] {SIDDS};
\addplot[draw=none, name path=lide_lb] table [x=correlation, y=F25] {\lide};
\addplot[draw=none, name path=lide_ub] table [x=correlation, y=F75] {\lide};
\addplot[red, opacity=0.3] fill between [of=lide_lb and lide_ub];

\pgfplotstableread{data/fig_correlated_vanderpol_crlb.dat}\crlb
\addplot[draw=none, domain = 1e-3:1, name path = crb] table[x=correlation, y expr=\thisrow{crlb_median}*1e-2] 
	{\crlb};
\addplot[draw=none, name path = crblb, domain=0:1] {1e-3}
	node [pos=0.5, anchor = south, rotate=0] {CRLB};
\addplot[gray, opacity =0.3] fill between [of=crblb and crb];

\end{groupplot}
\end{tikzpicture}

\caption[Correlated]{%
	SIDDS obtains the lower bound 
	when taking $\vma M = \vma \Sigma^{-1}$
	for correlated measurements described by~\cref{eq:correlated} with $\sigma=10^{-2}$.
	WLSOI refers to the weighted variant of LSOI described in~\cref{eq:weighted_lsoi}.
}
\label{fig:correlated}
\end{figure}

With SIDDS, 
we can also incorporate knowledge of noise correlation. 
Suppose that noise $\ve n_j$ at time $j$ is correlated between coordinates
with correlation $\rho\in [0, 1)$:
\begin{align}\label{eq:correlated}
	\ve n_j &\sim \mathcal{N}\left(
		\! \ve 0,
		\sigma^2\! \begin{bmatrix}
			1 & \rho \\ \rho & 1
		\end{bmatrix}
		\right)\! ,
	\ \text{if} \ \ve n_j \in \R^2; &
	\ve n_j &\sim \mathcal{N}\left(
		\! \ve 0,
		\sigma^2 \! \begin{bmatrix}
			1 & \rho & 0\\ \rho & 1 & 0 \\ 0 & 0 & 1
		\end{bmatrix}
		\right)\! ,
	\ \text{if} \ \ve n_j \in \R^3.
\end{align}
Then $\vve n \sim\mathcal{N}(\vve 0, \vma \Sigma)$ where
$\vma \Sigma$ is a block diagonal matrix consisting of $m$ repetitions of the block above.
Taking the weight $\vma M = \vma \Sigma^{-1}$ in SIDDS, 
we obtain near optimal estimates as shown in \cref{fig:correlated}.

Could we incorporate knowledge that noise has 
the distribution $\vve n \in \mathcal{N}(\vve 0, \vma \Sigma)$
into LSOI?
We are not aware of any existing work that does,
but we can using a weighted LSOI;
i.e., for some $\vma \Gamma$, solving
\begin{equation}\label{eq:weighted_lsoi}
	\min_{\vve c} \| \vma \Gamma [ \vma D \vve y - \vma \Phi(\vve y) \vve c]\|_2^2.
\end{equation}
In a linear estimation problem
$\min_{\vve c} \|\vma \Gamma[\vve y- \vma A \vve c]\|_2^2$
we choose $\vma \Gamma = \vma \Sigma^{-1/2}$
so that the residual is i.i.d.\ normally distributed 
(this is sometimes called \emph{whitening}).
Similarly, for LSOI linearizing around the true solution yields
$\vma \Gamma = \vma \Sigma^{-1/2}[ \vma D - \nabla \vma \Phi(\vve x) \tenmult2 \vve c^\star]^+$.
As seen in \cref{fig:correlated}
this approach can sometimes yield improvements over unweighted LSOI,
but the ill-conditioning of $\vma \Gamma$ 
limits the utility of this approach.
If we use a rank-truncated pseudoinverse
we can recover estimates that nearly obtain the CRLB
at the cost of an $\order( (md)^3)$ operation SVD;
this infeasible for large problems.
Ill-conditioning can be avoided by reformulating weighted LSOI  
as a constrained optimization problem,
\begin{equation}
	\begin{split} 
		\min_{\vve c, \vve r} \ & \| \vma \Sigma^{-1/2} \vve r\|_2  \\
		\text{s.t.} \ & [\vma D - \nabla \vma \Phi(\vve x)\tenmult2 \vve c^\star]\vve r =
			\vma D \vve y - \vma \Phi(\vve y) \vve c;
	\end{split}
\end{equation}
however, at this point we have reinvented SIDDS,
albeit using a different set of variables 
and a linearized constraint around the true values $\vve c^\star$ and $\vve x$.

\subsection{Large Data}
\begin{figure}

\centering
\begin{tikzpicture}
\begin{axis}[
	xmode = log,
	ymode = log,
	xlabel = data length $m$,
	ylabel = $\| \vve c - \vve c^\star \|_2$,
	width = 0.97\linewidth,
	height = 0.4\linewidth,
	ytickten = {-7,-6, ..., 5},
	ymajorgrids,
	ymin = 1e-3,
	ymax = 100,
	xmin = 0,
	xmax = 100000,
	title = {Recovering the Lorenz 63 Attractor with Increasing Data},
]
\pgfplotstableread{data/fig_length_lorenz63_lsopinf.dat}\lsopinf
\addplot[blue] table [x=length, y=F50] {\lsopinf}
	node [pos=0.63, anchor = south, rotate=0, yshift=2pt] {LSOI};
\addplot[draw=none, name path=lsopinf_lb] table [x=length, y=F25] {\lsopinf};
\addplot[draw=none, name path=lsopinf_ub] table [x=length, y=F75] {\lsopinf};
\addplot[blue, opacity=0.3] fill between [of=lsopinf_lb and lsopinf_ub];

\pgfplotstableread{data/fig_length_lorenz63_lide.dat}\lide
\addplot[red] table [x=length, y=F50] {\lide}
	node [pos=.6, anchor = south, rotate=0, yshift=3pt] {SIDDS};
\addplot[draw=none, name path=lide_lb] table [x=length, y=F25] {\lide};
\addplot[draw=none, name path=lide_ub] table [x=length, y=F75] {\lide};
\addplot[red, opacity=0.3] fill between [of=lide_lb and lide_ub];

\pgfplotstableread{data/fig_length_crlb.dat}\crlb
\addplot[draw=none, domain = 1e-3:1, name path = crb] table[x=m, y expr = 1e-2 * \thisrow{crlb_med}] {\crlb};
\addplot[draw=none, name path = crblb, domain=1e2:21544] {1e-3}
	node [pos=0.0, anchor = south west, rotate=0] {CRLB};
\addplot[gray, opacity =0.3] fill between [of=crblb and crb];

\end{axis}
\end{tikzpicture}
\caption[Data]{
SIDDS enables the accurate recovery of dynamical systems 
with less data than LSOI.
This example considers the Lorenz 63 attractor with sample rate $\delta = 10^{-2}$
and additive i.i.d.\ normally distributed noise with standard deviation $\sigma = 10^{-2}$,
using $\vma M = \ma I$ in SIDDS.
Numerics limit the ability to compute the CRLB for large $m$.
}
\label{fig:length}
\end{figure}
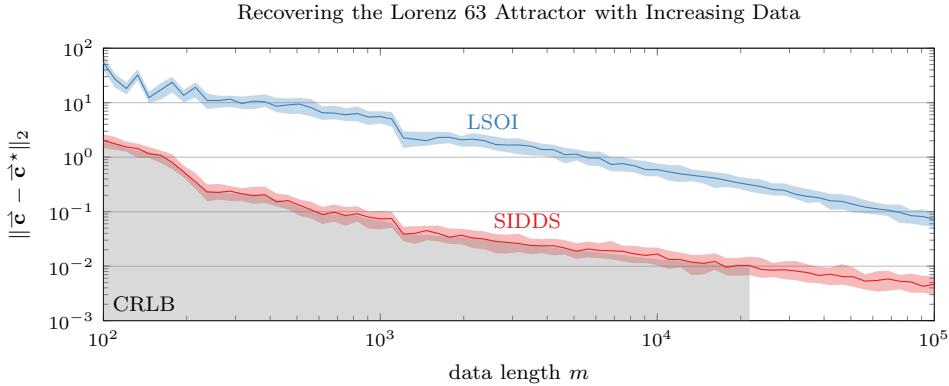

The previous example showed that SIDDS recovers
more accurate estimates than LSOI.
Another interpretation of this result is that SIDDS 
obtains similarly accurate estimates
using less data.
\Cref{fig:length} illustrates this point
by recovering the chaotic Lorenz 63 attractor
with increasing amounts of data.
We observe SIDDS obtains roughly the same accuracy
as LSOI using ten times less data.
This example also serves as a stress-test of SIDDS,
illustrating that the algorithm scales to large scale problems
while still approximately obtaining the CRLB.


\subsection{Comparing SIDDS to Other Algorithms}
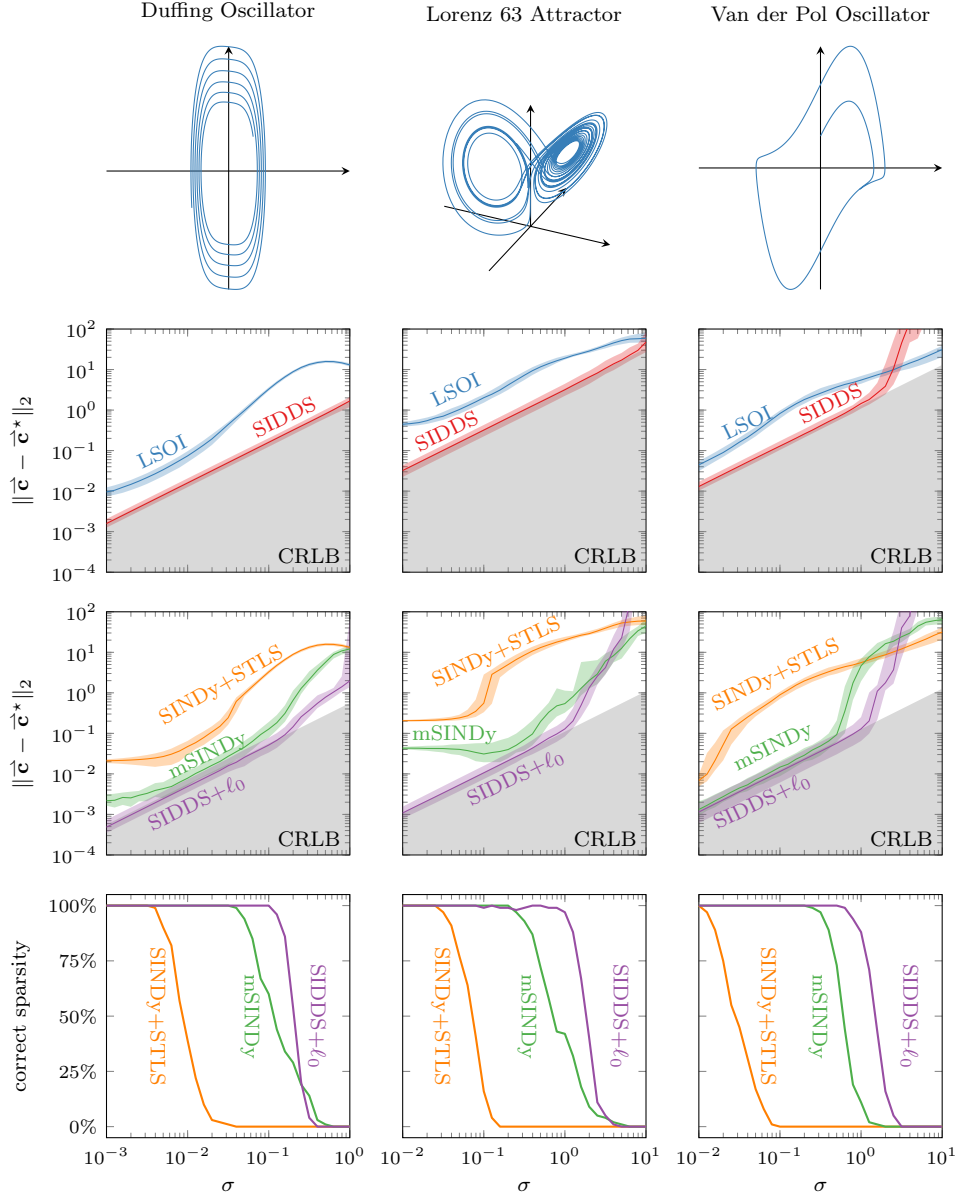
\begin{figure}
\centering

\begin{tikzpicture}
\begin{groupplot}[
	group style = {group size = 3 by 4,
		horizontal sep = 2em,
		vertical sep = 1.5em,
		},
	width = 0.37\linewidth,
	height = 0.37\linewidth,
]

\nextgroupplot[
	axis equal,
	axis lines= middle,
	xtick = \empty, 
	ytick = \empty,
	title = Duffing Oscillator,
	clip = false,
]
	\addplot [smooth, blue] table [x=x, y=y] {data/fig_noise_duffing_state.dat};

\nextgroupplot[
	axis equal,
	axis lines= middle,
	xtick = \empty, 
	ytick = \empty,
	ztick = \empty,
	title = Lorenz 63 Attractor,
	clip = false,
]
	\addplot3 [smooth, blue] table [x=x, y=y, z=z] {data/fig_noise_lorenz63_state.dat};

\nextgroupplot[
	axis equal,
	axis lines= middle,
	xtick = \empty, 
	ytick = \empty,
	ztick = \empty,
	title = Van der Pol Oscillator,
	clip = false,
]
	\addplot [smooth, blue] table [x=x, y=y] {data/fig_noise_vanderpol_state.dat};


\nextgroupplot[
	ylabel = $\| \vve c - \vve c^\star\|_2$,
	ytickten = {-7,-6, ..., 5},
	ymin = 1e-4,
	ymax = 1e2,
	xmode = log,
	ymode = log,
	xticklabels = {,,},
	]

\pgfplotstableread{data/fig_noise_duffing_lsopinf.dat}\lsopinf
\addplot[blue] table [x=noise, y=F50] {\lsopinf}
	node [pos=0.1, anchor = south west, rotate=26, yshift=0pt] {LSOI};
\addplot[draw=none, name path=lsopinf_lb] table [x=noise, y=F25] {\lsopinf};
\addplot[draw=none, name path=lsopinf_ub] table [x=noise, y=F75] {\lsopinf};
\addplot[blue, opacity=0.3] fill between [of=lsopinf_lb and lsopinf_ub];

\pgfplotstableread{data/fig_noise_duffing_lide_smooth.dat}\lide
\addplot[red] table [x=noise, y=F50] {\lide}
	node [pos=0.6, anchor = south west, rotate=26] {SIDDS};
\addplot[draw=none, name path=lide_lb] table [x=noise, y=F25] {\lide};
\addplot[draw=none, name path=lide_ub] table [x=noise, y=F75] {\lide};
\addplot[red, opacity=0.3] fill between [of=lide_lb and lide_ub];

\addplot[draw=none, domain = 1e-3:1, name path = crb] {x*1.68388448};
\addplot[draw=none, name path = crblb, domain=1e-3:1] {1e-4}
	node [pos=1, anchor = south east, rotate=0] {CRLB};
\addplot[gray, opacity =0.3] fill between [of=crblb and crb];

\nextgroupplot[
	ytickten = {-7,-6, ..., 5},
	ymin = 1e-4,
	ymax = 1e2,
	xmode = log,
	ymode = log,
	xticklabels = {,,},
	yticklabels = {,,},
	]

\pgfplotstableread{data/fig_noise_lorenz63_lsopinf.dat}\lsopinf
\addplot[blue] table [x=noise, y=F50] {\lsopinf}
	node [pos=0.1, anchor = south west, rotate=26, yshift=0pt] {LSOI};
\addplot[draw=none, name path=lsopinf_lb] table [x=noise, y=F25] {\lsopinf};
\addplot[draw=none, name path=lsopinf_ub] table [x=noise, y=F75] {\lsopinf};
\addplot[blue, opacity=0.3] fill between [of=lsopinf_lb and lsopinf_ub];

\pgfplotstableread{data/fig_noise_lorenz63_lide_smooth.dat}\lide
\addplot[red] table [x=noise, y=F50] {\lide}
	node [pos=0.05, anchor = south west, rotate=26] {SIDDS};
\addplot[draw=none, name path=lide_lb] table [x=noise, y=F25] {\lide};
\addplot[draw=none, name path=lide_ub] table [x=noise, y=F75] {\lide};
\addplot[red, opacity=0.3] fill between [of=lide_lb and lide_ub];

\addplot[draw=none, domain = 1e-2:1e1, name path = crb] {x*3.35201717};
\addplot[draw=none, name path = crblb, domain=1e-2:1e1] {1e-4}
	node [pos=1, anchor = south east, rotate=0] {CRLB};
\addplot[gray, opacity =0.3] fill between [of=crblb and crb];

\nextgroupplot[
	ytickten = {-7,-6, ..., 5},
	ymin = 1e-4,
	ymax = 1e2,
	xmode = log,
	ymode = log,
	xticklabels = {,,},
	yticklabels = {,,},
	]

\pgfplotstableread{data/fig_noise_vanderpol_lsopinf.dat}\lsopinf
\addplot[blue] table [x=noise, y=F50] {\lsopinf}
	node [pos=0.1, anchor = south west, rotate=26, yshift=0pt] {LSOI};
\addplot[draw=none, name path=lsopinf_lb] table [x=noise, y=F25] {\lsopinf};
\addplot[draw=none, name path=lsopinf_ub] table [x=noise, y=F75] {\lsopinf};
\addplot[blue, opacity=0.3] fill between [of=lsopinf_lb and lsopinf_ub];

\pgfplotstableread{data/fig_noise_vanderpol_lide_smooth.dat}\lide
\addplot[red] table [x=noise, y=F50] {\lide}
	node [pos=0.3, anchor = south west, rotate=26] {SIDDS};
\addplot[draw=none, name path=lide_lb] table [x=noise, y=F25] {\lide};
\addplot[draw=none, name path=lide_ub] table [x=noise, y=F75] {\lide};
\addplot[red, opacity=0.3] fill between [of=lide_lb and lide_ub];

\addplot[draw=none, domain = 1e-2:1e1, name path = crb] {x*1.30300513};
\addplot[draw=none, name path = crblb, domain=1e-2:1e1] {1e-4}
	node [pos=1, anchor = south east, rotate=0] {CRLB};
\addplot[gray, opacity =0.3] fill between [of=crblb and crb];


\nextgroupplot[
	ytickten = {-7,-6, ..., 5},
	ymin = 1e-4,
	ymax = 1e2,
	xmode = log,
	ymode = log,
	ylabel = $\| \vve c - \vve c^\star \|_2$,
	xticklabels = {,,},
	]

\pgfplotstableread{data/fig_noise_duffing_sindy.dat}\sindy
\addplot[orange] table [x=noise, y=F50] {\sindy}
	node [pos=0.5, anchor = south, rotate=30, yshift=2pt] {SINDy+STLS};
\addplot[draw=none, name path=sindy_lb] table [x=noise, y=F25] {\sindy};
\addplot[draw=none, name path=sindy_ub] table [x=noise, y=F75] {\sindy};
\addplot[orange, opacity=0.3] fill between [of=sindy_lb and sindy_ub];

\pgfplotstableread{data/fig_noise_duffing_mod_sindy.dat}\msindy
\addplot[green] table [x=noise, y=F50] {\msindy}
	node [pos=0.5, anchor = south east, rotate=26, yshift=-5pt] {mSINDy};
\addplot[draw=none, name path=msindy_lb] table [x=noise, y=F25] {\msindy};
\addplot[draw=none, name path=msindy_ub] table [x=noise, y=F75] {\msindy};
\addplot[green, opacity=0.3] fill between [of=msindy_lb and msindy_ub];

\pgfplotstableread{data/fig_noise_duffing_lide_irls0_smooth_polish_alpha_0.01.dat}\lidezerof
\addplot[purple] table [x=noise, y=F50] {\lidezerof}
	node [pos=.1, anchor = north west, rotate=26] {SIDDS+$\ell_0$};
\addplot[draw=none, name path=lidezerof_lb] table [x=noise, y=F25] {\lidezerof};
\addplot[draw=none, name path=lidezerof_ub] table [x=noise, y=F75] {\lidezerof};
\addplot[purple, opacity=0.3] fill between [of=lidezerof_lb and lidezerof_ub];

\addplot[draw=none, domain = 1e-3:1e0, name path = crb] {x*0.56065196};
\addplot[draw=none, name path = crblb, domain=1e-3:1e0] {1e-4}
	node [pos=1, anchor = south east, rotate=0] {CRLB};
\addplot[gray, opacity =0.3] fill between [of=crblb and crb];

\nextgroupplot[
	ytickten = {-7,-6, ..., 5},
	ymin = 1e-4,
	ymax = 1e2,
	xmode = log,
	ymode = log,
	xticklabels = {,,},
	yticklabels = {,,},
	]

\pgfplotstableread{data/fig_noise_lorenz63_sindy.dat}\sindy
\addplot[orange] table [x=noise, y=F50] {\sindy}
	node [pos=0.5, anchor = south, rotate=25, xshift=0pt, yshift=-1pt] {SINDy+STLS};
\addplot[draw=none, name path=sindy_lb] table [x=noise, y=F25] {\sindy};
\addplot[draw=none, name path=sindy_ub] table [x=noise, y=F75] {\sindy};
\addplot[orange, opacity=0.3] fill between [of=sindy_lb and sindy_ub];

\pgfplotstableread{data/fig_noise_lorenz63_mod_sindy.dat}\msindy
\addplot[green] table [x=noise, y=F50] {\msindy}
	node [pos=0.0, anchor = south west, rotate=0, xshift=0pt, yshift=-2pt] {mSINDy};
\addplot[draw=none, name path=msindy_lb] table [x=noise, y=F25] {\msindy};
\addplot[draw=none, name path=msindy_ub] table [x=noise, y=F75] {\msindy};
\addplot[green, opacity=0.3] fill between [of=msindy_lb and msindy_ub];

\pgfplotstableread{data/fig_noise_lorenz63_lide_irls0_smooth_polish_alpha_0.5.dat}\lidezero
\addplot[purple] table [x=noise, y=F50] {\lidezero}
	node [pos=.1, anchor = north west, rotate=26] {SIDDS+$\ell_0$};
\addplot[draw=none, name path=lidezero_lb] table [x=noise, y=F25] {\lidezero};
\addplot[draw=none, name path=lidezero_ub] table [x=noise, y=F75] {\lidezero};
\addplot[purple, opacity=0.3] fill between [of=lidezero_lb and lidezero_ub];

\addplot[draw=none, domain = 1e-2:1e1, name path = crb] {x*0.11003616};
\addplot[draw=none, name path = crblb, domain=1e-2:1e1] {1e-4}
	node [pos=1, anchor = south east, rotate=0] {CRLB};
\addplot[gray, opacity =0.3] fill between [of=crblb and crb];


\nextgroupplot[
	ytickten = {-7,-6, ..., 5},
	ymin = 1e-4,
	ymax = 1e2,
	xmode = log,
	ymode = log,
	xticklabels = {,,},
	yticklabels = {,,},
	]

\pgfplotstableread{data/fig_noise_vanderpol_sindy.dat}\sindy
\addplot[orange] table [x=noise, y=F50] {\sindy}
	node [pos=0.5, anchor = south, rotate=26, yshift=0pt] {SINDy+STLS};
\addplot[draw=none, name path=sindy_lb] table [x=noise, y=F25] {\sindy};
\addplot[draw=none, name path=sindy_ub] table [x=noise, y=F75] {\sindy};
\addplot[orange, opacity=0.3] fill between [of=sindy_lb and sindy_ub];

\pgfplotstableread{data/fig_noise_vanderpol_mod_sindy.dat}\msindy
\addplot[green] table [x=noise, y=F50] {\msindy}
	node [pos=0.4, anchor = south east, rotate=26, yshift=0pt] {mSINDy};
\addplot[draw=none, name path=msindy_lb] table [x=noise, y=F25] {\msindy};
\addplot[draw=none, name path=msindy_ub] table [x=noise, y=F75] {\msindy};
\addplot[green, opacity=0.3] fill between [of=msindy_lb and msindy_ub];

\pgfplotstableread{data/fig_noise_vanderpol_lide_irls0_smooth_polish_alpha_1.0.dat}\lidezero
\addplot[purple] table [x=noise, y=F50] {\lidezero}
	node [pos=0, anchor = north west, rotate=26] {SIDDS+$\ell_0$};
\addplot[draw=none, name path=lidezero_lb] table [x=noise, y=F25] {\lidezero};
\addplot[draw=none, name path=lidezero_ub] table [x=noise, y=F75] {\lidezero};
\addplot[purple, opacity=0.3] fill between [of=lidezero_lb and lidezero_ub];

\addplot[draw=none, domain = 1e-2:1e1, name path = crb] {x*0.12670384};
\addplot[draw=none, name path = crblb, domain=1e-2:1e1] {1e-4}
	node [pos=1, anchor = south east, rotate=0] {CRLB};
\addplot[gray, opacity =0.3] fill between [of=crblb and crb];


\nextgroupplot[
	xmode = log,
	ymin = -.05, ymax=1.05,
	ylabel = correct sparsity,
	yticklabel={\pgfmathparse{int(100*\tick)} \pgfmathresult\%},
	ytick = {0,0.25, 0.5, 0.75, 1},
	xlabel=$\sigma$,
	]

\pgfplotstableread{data/fig_noise_duffing_sindy.dat}\sindy
\addplot[orange, thick] table [x=noise, y=exact_sparse] {\sindy};
\node at (axis cs:0.8e-2,0.5)      [anchor=north, rotate=-90, orange, yshift=-2pt] {SINDy+STLS};

\pgfplotstableread{data/fig_noise_duffing_mod_sindy.dat}\msindy
\addplot[green, thick] table [x=noise, y=exact_sparse] {\msindy};
\node at (axis cs:1e-1,0.5)      [anchor=north, rotate=-90, green] {mSINDy};

\pgfplotstableread{data/fig_noise_duffing_lide_irls0_smooth_polish_alpha_0.01.dat}\lidezerof
\addplot[purple, thick] table [x=noise, y=exact_sparse] {\lidezerof};
\node at (axis cs:5e-1,0.5)      [anchor=south, rotate=-90, purple, yshift=-10pt] {SIDDS+$\ell_0$};

\nextgroupplot[
	xmode = log,
	ymin = -.05, ymax=1.05,
	yticklabel={\pgfmathparse{int(100*\tick)} \pgfmathresult\%},
	ytick = {0,0.25, 0.5, 0.75, 1},
	xlabel=$\sigma$,
	yticklabels = {,,},
	]

\pgfplotstableread{data/fig_noise_lorenz63_sindy.dat}\sindy
\addplot[orange, thick] table [x=noise, y=exact_sparse] {\sindy};
\node at (axis cs:5e-2,0.5)      [anchor=north, rotate=-90, orange] {SINDy+STLS};

\pgfplotstableread{data/fig_noise_lorenz63_mod_sindy.dat}\msindy
\addplot[green, thick] table [x=noise, y=exact_sparse] {\msindy};
\node at (axis cs:5e-1,0.5)      [anchor=north, rotate=-90, green, yshift=2pt] {mSINDy};

\pgfplotstableread{data/fig_noise_lorenz63_lide_irls0_smooth_polish_alpha_0.5.dat}\lidezero
\addplot[purple, thick] table [x=noise, y=exact_sparse] {\lidezero};
\node at (axis cs:6e0,0.5)      [anchor=south, rotate=-90, purple, yshift=-10pt] {SIDDS+$\ell_0$};

\nextgroupplot[
	xmode = log,
	ymin = -.05, ymax=1.05,
	yticklabel={\pgfmathparse{int(100*\tick)} \pgfmathresult\%},
	ytick = {0,0.25, 0.5, 0.75, 1},
	xlabel=$\sigma$,
	yticklabels = {,,},
	]
\pgfplotstableread{data/fig_noise_vanderpol_sindy.dat}\sindy
\addplot[orange, thick] table [x=noise, y=exact_sparse] {\sindy};
\node at (axis cs:1e-1,0.5)      [anchor=north, rotate=-90, orange, yshift=2pt] {SINDy+STLS};

\pgfplotstableread{data/fig_noise_vanderpol_mod_sindy.dat}\msindy
\addplot[green, thick] table [x=noise, y=exact_sparse] {\msindy};
\node at (axis cs:5e-1,0.5)  [anchor=north, rotate=-90, green] {mSINDy};

\pgfplotstableread{data/fig_noise_vanderpol_lide_irls0_smooth_polish_alpha_1.0.dat}\lidezero
\addplot[purple, thick] table [x=noise, y=exact_sparse] {\lidezero};
\node at (axis cs:5e0,0.5)   [anchor=south, rotate=-90, purple, yshift=-10pt] {SIDDS+$\ell_0$};

\end{groupplot}
\end{tikzpicture}
\caption[Noise]{%
	Both the SIDDS and the sparsity promoting SIDDS+$\ell_0$
	obtain better estimates than competing methods.
	For each test problem, we consider how each method
	performs as the standard deviation $\sigma$ of the 
	additive i.i.d.\ normally distributed noise increases.
	%
	The second row considers unregularized algorithms: SIDDS (with $\vma M = \vma I$) and LSOI.
	The third row considers sparsity promoting algorithms:
	SINDy using STLS as implemented in PySINDy~\cite{pysindy20,pysindy21} (SINDy+STLS),
	Modified SINDy as implemented by~\cite{KBK20x}  (mSINDy),
	and SIDDS+$\ell_0$ (with $\vma M = \vma I$).
	The CRLB for this row is computed assuming the correct sparsity structure in $\vve c$
	yielding a smaller lower bound than the preceding row as discussed in \cref{sec:crb:sparse}.
	%
	%
	%
	Here we use sample rate $\delta=10^{-2}$ 
	and $m$ of 1000, 2000, and 1000
	for Duffing, Lorenz 63, and Van der Pol respectively.
	For SINDy+STLS we use the default parameters in PySINDy.
	For mSINDy we use parameters from~\cite{KBK20x}:
	truncation parameter $\lambda$ of $0.05$, $0.2$, and $0.1$
	and ADAM steps per outer iteration of 5000, 15000, and 5000 respectively.
	For SIDDS+$\ell_0$, we choose the penalty $\alpha$ of $0.01$, $0.5$, and $1$ respectively
	chosen via a coarse optimization to maximize the probability 
	of recovering the correct sparsity structure
	with large noise.
}
\label{fig:noise}
\end{figure}

As a final example,
\cref{fig:noise} compares algorithms 
for estimating dynamical systems
perform with increasing levels of noise.
We separate these algorithms into two classes:
those without a sparsity promoting constraint (LSOI and SIDDS)
and those with a sparsity promoting constraint
(SINDy+STLS~\cite{BPK16,pysindy20}, Modified SINDy~\cite{KBK20x},
and SIDDS+$\ell_0$)
since sparsity promotion allows for an improved recovery
and lowers the CRLB as discussed in \cref{sec:crb:sparse}.
These experiments illustrate several key points.
First, SIDDS and SIDDS+$\ell_0$ approximately obtain the Cram\'er-Rao lower bound for small noise $\sigma$;
for large noise, these algorithms identify a local minimizer far away the true solution
leading the estimates to detach from the lower bound.
Second, while the denoising penalty introduced by mSINDy 
does improve estimates compared to SINDy+STLS,
it does not obtain the lower bound.
Third, SIDDS+$\ell_0$ is able to correctly identify the sparsity structure
for larger noise than SINDy+STLS and mSINDy.

\section{Discussion\label{sec:discussion}}
Here we have shown how to practically identify and denoise a dynamical system
using SIDDS and how to incorporate sparsity promotion in SIDDS+$\ell_0$.
This algorithm yields estimates obtaining the Cram\'er-Rao lower bound
for small noise, outperforming existing algorithms.
We anticipate there many possible avenues for improvement and extension of SIDDS;
e.g., better initialization through smoothing techniques described in~\cite{CPD22x}
and incorporating nonautonomous input as in~\cite{KKB18}.



\bibliographystyle{siamplain}
\bibliography{abbrevjournals,master}

\begin{thebibliography}{10}

\bibitem{AVBG17}
{\sc N.~Alger, U.~Villa, T.~Bui-Thanh, and O.~Ghattas}, {\em A data scalable
  augmented {Lagrangian} {KKT} preconditioner for large-scale inverse
  problems}, SIAM J. Sci. Comput., 39 (2017), pp.~A2365--A2393,
  \url{https://doi.org/10.1137/16m1084365}.

\bibitem{ABN16}
{\sc M.~Asch, M.~Bocquet, and M.~Nodet}, {\em Data Assimilation}, Society for
  Industrial and Applied Mathematics, Dec. 2016,
  \url{https://doi.org/10.1137/1.9781611974546}.

\bibitem{BE09}
{\sc Z.~Ben-Haim and Y.~C. Eldar}, {\em On the constrained {Cram\'er--Rao}
  bound with a singular {Fisher} information matrix}, IEEE Signal Process.
  Lett., 16 (2009), pp.~453--456,
  \url{https://doi.org/10.1109/lsp.2009.2016831}.

\bibitem{Bjo96}
{\sc {\AA}.~Bj{\"o}rck}, {\em Numerical Methods for Least Squares Problems},
  SIAM, Philadelphia, PA, 1996.

\bibitem{BD09}
{\sc T.~Blumensath and M.~E. Davies}, {\em Iterative hard thresholding for
  compressed sensing}, Appl. Comput. Harmon. A., 27 (2009), pp.~265--274,
  \url{https://doi.org/10.1016/j.acha.2009.04.002}.

\bibitem{BPK16}
{\sc S.~L. Brunton, J.~L. Proctor, and J.~N. Kutz}, {\em Discovering governing
  equations from data by sparse identification of nonlinear dynamical systems},
  Proc. Natl. Acad. Sci. USA, 113 (2016), pp.~3932--3937.

\bibitem{CWB08}
{\sc E.~J. Cand{\`e}s, M.~B. Wakin, and S.~P. Boyd}, {\em Enforcing sparsity by
  reweighted $\ell_1$ minimization}, J. Fourier Anal. Appl., 14 (2008),
  pp.~877--905, \url{https://doi.org/10.1007/s00041-008-9045-x}.

\bibitem{CY08}
{\sc R.~Chartrand and W.~Yin}, {\em Iteratively reweighted algorithms for
  compressive sensing}, in 2008 IEEE International Conference on Acoustics,
  Speech and Signal Processing, IEEE, Mar. 2008,
  \url{https://doi.org/10.1109/icassp.2008.4518498}.

\bibitem{CPD21}
{\sc A.~Cortiella, K.-C. Park, and A.~Doostan}, {\em Sparse identification of
  nonlinear dynamical systems via reweighted {$\ell_1$}-regularized least
  squares}, Comput. Method. Appl. M., 376 (2021), p.~113620,
  \url{https://doi.org/10.1016/j.cma.2020.113620}.

\bibitem{CPD22x}
{\sc A.~Cortiella, K.-C. Park, and A.~Doostan}, {\em A priori denoising
  strategies for sparse identification of nonlinear dynamical systems: A
  comparative study}, 2022, \url{https://arxiv.org/abs/arxiv:2201.12683v1}.

\bibitem{DDFG10}
{\sc I.~Daubechies, R.~DeVore, M.~Fornasier, and C.~S. Güntürk}, {\em
  Iteratively reweighted least squares minimization for sparse recovery}, Comm.
  Pure Appl. Math., 63 (2010), pp.~1--38,
  \url{https://doi.org/10.1002/cpa.20303}.

\bibitem{pysindy20}
{\sc B.~de~Silva, K.~Champion, M.~Quade, J.-C. Loiseau, J.~Kutz, and
  S.~Brunton}, {\em {PySINDy:} {A} {P}ython package for the sparse
  identification of nonlinear dynamical systems from data}, JOSS, 5 (2020),
  p.~2104, \url{https://doi.org/10.21105/joss.02104}.

\bibitem{superlu99}
{\sc J.~W. Demmel, S.~C. Eisenstat, J.~R. Gilbert, X.~S. Li, and J.~W.~H. Liu},
  {\em A supernodal approach to sparse partial pivoting}, SIAM J. Matrix Anal.
  \& Appl., 20 (1999), pp.~720--755,
  \url{https://doi.org/10.1137/s0895479895291765}.

\bibitem{GP73}
{\sc G.~H. Golub and V.~Pereyra}, {\em The differentiation of pseudo-inverses
  and nonlinear least squares problems whose variables separate}, SIAM J.
  Numer. Anal., 10 (1973), pp.~413--432, \url{https://doi.org/10.1137/0710036}.

\bibitem{GL13}
{\sc G.~H. Golub and C.~F. Van~Loan}, {\em Matrix Computations}, Johns Hopkins
  University Press, Baltimore, MD, fourth~ed., 2013.

\bibitem{HA01}
{\sc E.~Haber and U.~M. Ascher}, {\em Preconditioned all-at-once methods for
  large, sparse parameter estimation problems}, Inverse Probl., 17 (2001),
  pp.~1847--1864, \url{https://doi.org/10.1088/0266-5611/17/6/319}.

\bibitem{EAO00}
{\sc E.~Haber, U.~M. Ascher, and D.~Oldenburg}, {\em On optimization techniques
  for solving nonlinear inverse problems}, Inverse Probl., 16 (2000),
  pp.~1263--1280, \url{https://doi.org/10.1088/0266-5611/16/5/309}.

\bibitem{Hei13x}
{\sc M.~Heinkenschloss}, {\em Numerical solution of implicitly constrained
  optimization problems}, Tech. Report TR08-05, Rice University, 2013,
  \url{https://hdl.handle.net/1911/102087}.

\bibitem{HJ12}
{\sc R.~A. Horn and C.~R. Johnson}, {\em Matrix Analysis}, Cambridge University
  Press, second~ed., 2012, \url{https://doi.org/10.1017/cbo9781139020411}.

\bibitem{KBK20x}
{\sc K.~Kaheman, S.~L. Brunton, and J.~N. Kutz}, {\em Automatic differentiation
  to simultaneously identify nonlinear dynamics and extract noise probability
  distributions from data}, 2020, \url{https://arxiv.org/abs/2009.08810v2}.

\bibitem{KKB18}
{\sc E.~Kaiser, J.~Kutz, and S.~Brunton}, {\em Sparse identification of
  nonlinear dynamics for model predictive control in the low-data limit}, Proc.
  R. Soc. A., 474 (2018), p.~20180335,
  \url{https://doi.org/10.1098/rspa.2018.0335},
  \url{https://doi.org/10.1098/rspa.2018.0335}.

\bibitem{pysindy21}
{\sc A.~A. Kaptanoglu, B.~M. de~Silva, U.~Fasel, K.~Kaheman, A.~J. Goldschmidt,
  J.~L. Callaham, C.~B. Delahunt, Z.~G. Nicolaou, K.~Champion, J.-C. Loiseau,
  J.~N. Kutz, and S.~L. Brunton}, {\em {PySINDy}: {A} comprehensive {P}ython
  package for robust sparse system identification}, arXiv preprint
  arXiv:2111.08481,  (2021).

\bibitem{KB09}
{\sc T.~G. Kolda and B.~W. Bader}, {\em Tensor decompositions and
  applications}, SIAM Rev., 51 (2009), pp.~455--500,
  \url{https://doi.org/10.1137/07070111X}.

\bibitem{LXY13}
{\sc M.-J. Lai, Y.~Xu, and W.~Yin}, {\em Improved iteratively reweighted least
  squares for unconstrained smoothed $\ell_q$ minimization}, SIAM J. Numer.
  Anal., 51 (2013), pp.~927--957, \url{https://doi.org/10.1137/110840364}.

\bibitem{LY11}
{\sc X.~Liu and Y.~Yuan}, {\em A sequential quadratic programming method
  without a penalty function or a filter for nonlinear equality constrained
  optimization}, SIAM J. Optim., 21 (2011), pp.~545--571,
  \url{https://doi.org/10.1137/080739884}.

\bibitem{PW16}
{\sc B.~Peherstofer and K.~Willcox}, {\em Data-driven operator inference for
  nonintrusive projection-based model reduction}, Comput. Methods Appl. Mech.
  Engrg., 306 (2016), pp.~196--215,
  \url{https://doi.org/10.1016/j.cma.2016.03.025}.

\bibitem{RKB19}
{\sc S.~H. Rudy, J.~N. Kutz, and S.~L. Brunton}, {\em Deep learning of dynamics
  and signal-noise decomposition with time-stepping constraints}, J. Comput.
  Phys., 396 (2019), pp.~483--506,
  \url{https://doi.org/10.1016/j.jcp.2019.06.056}.

\bibitem{SW89}
{\sc G.~A.~F. Seber and C.~J. Wild}, {\em Nonlinear Regression},
  Wiley-Interscience, 1989.

\bibitem{SLS21}
{\sc F.~Sun, Y.~Liu, and H.~Sun}, {\em Physics-informed spline learning for
  nonlinear dynamics discovery}, in Proceedings of the Thirtieth International
  Joint Conference on Artificial Intelligence, International Joint Conferences
  on Artificial Intelligence Organization, Aug. 2021,
  \url{https://doi.org/10.24963/ijcai.2021/283}.

\bibitem{TW17}
{\sc G.~Tran and R.~Ward}, {\em Exact recovery of chaotic systems from highly
  corrupted data}, Multiscale Model. Simul., 15 (2017), pp.~1108--1129,
  \url{https://doi.org/10.1137/16m1086637}.

\bibitem{Wri05}
{\sc S.~J. Wright}, {\em An algorithm for degenerate nonlinear programming with
  rapid local convergence}, SIAM J. Optim., 15 (2005), pp.~673--696,
  \url{https://doi.org/10.1137/030601235}.

\bibitem{WX19}
{\sc K.~Wu and D.~Xiu}, {\em Numerical aspects for approximating governing
  equations using data}, J. Comput. Phys., 384 (2019), pp.~200--221,
  \url{https://doi.org/10.1016/j.jcp.2019.01.030},
  \url{https://doi.org/10.1016/j.jcp.2019.01.030}.

\end{thebibliography}

\end{document}